\documentclass[11pt,twoside]{atmp}

\usepackage{amsmath,amssymb}
\usepackage[all]{xy}
\usepackage{color}

\usepackage{amscd}
\usepackage{amssymb}
\usepackage{array}
\usepackage{color}
\usepackage{youngtab}

\newtheorem{thm}{Theorem}[section]
\newtheorem{theorem}[thm]{Theorem}

\newtheorem{proposition}[thm]{Proposition}

\theoremstyle{definition}
\newtheorem{definition}[thm]{Definition}
\newtheorem{example}[thm]{Example}

\newtheorem{observation}[thm]{Observation}

\newtheorem{remark}[thm]{Remark}

\begin{document}

\newcommand{\comment}[1]{{\color{blue}\rule[-0.5ex]{2pt}{2.5ex}}
\marginpar{\scriptsize\begin{flushleft}\color{blue}#1\end{flushleft}}}

\newcommand{\be}{\begin{equation}}
\newcommand{\ee}{\end{equation}}
\newcommand{\bea}{\begin{eqnarray}}
\newcommand{\eea}{\end{eqnarray}}
\newcommand{\bean}{\begin{eqnarray*}}
\newcommand{\eean}{\end{eqnarray*}}

\newcommand{\id}{\relax{\rm 1\kern-.28em 1}}
\newcommand{\R}{\mathbb{R}}
\newcommand{\C}{\mathbb{C}}
\newcommand{\Z}{\mathbb{Z}}
\newcommand{\g}{\mathfrak{G}}
\newcommand{\e}{\epsilon}

\newcommand{\hs}{\hfill\square}
\newcommand{\hbs}{\hfill\blacksquare}

\newcommand{\bp}{\mathbf{p}}
\newcommand{\bmax}{\mathbf{m}}
\newcommand{\bT}{\mathbf{T}}
\newcommand{\bU}{\mathbf{U}}
\newcommand{\bP}{\mathbf{P}}
\newcommand{\bA}{\mathbf{A}}
\newcommand{\bm}{\mathbf{m}}
\newcommand{\bIP}{\mathbf{I_P}}

\newcommand{\cA}{\mathcal{A}}
\newcommand{\cB}{\mathcal{B}}
\newcommand{\cC}{\mathcal{C}}
\newcommand{\cI}{\mathcal{I}}
\newcommand{\cO}{\mathcal{O}}
\newcommand{\cG}{\mathcal{G}}
\newcommand{\cJ}{\mathcal{J}}
\newcommand{\cF}{\mathcal{F}}
\newcommand{\cP}{\mathcal{P}}
\newcommand{\ep}{\mathcal{E}}
\newcommand{\E}{\mathcal{E}}
\newcommand{\cH}{\mathcal{O}}
\newcommand{\cPO}{\mathcal{PO}}
\newcommand{\cl}{\ell}
\newcommand{\cFG}{\mathcal{F}_{\mathrm{G}}}
\newcommand{\cHG}{\mathcal{H}_{\mathrm{G}}}
\newcommand{\Gal}{G_{\mathrm{al}}}
\newcommand{\cQ}{G_{\mathcal{Q}}}

\newcommand{\ri}{\mathrm{i}}
\newcommand{\re}{\mathrm{e}}
\newcommand{\rd}{\mathrm{d}}
\newcommand{\rSt}{\mathrm{St}}
\newcommand{\rGL}{\mathrm{GL}}
\newcommand{\rSU}{\mathrm{SU}}
\newcommand{\rSL}{\mathrm{SL}}
\newcommand{\rSO}{\mathrm{SO}}
\newcommand{\rOSp}{\mathrm{OSp}}
\newcommand{\rSpin}{\mathrm{Spin}}
\newcommand{\rsl}{\mathrm{sl}}
\newcommand{\rM}{\mathrm{M}}
\newcommand{\rdiag}{\mathrm{diag}}
\newcommand{\rP}{\mathrm{P}}
\newcommand{\rdeg}{\mathrm{deg}}
\newcommand{\rStab}{\mathrm{Stab}}

\newcommand{\M}{\mathrm{M}}
\newcommand{\End}{\mathrm{End}}
\newcommand{\Hom}{\mathrm{Hom}}
\newcommand{\diag}{\mathrm{diag}}
\newcommand{\rspan}{\mathrm{span}}
\newcommand{\rank}{\mathrm{rank}}
\newcommand{\Gr}{\mathrm{Gr}}
\newcommand{\ber}{\mathrm{Ber}}

\newcommand{\fsl}{\mathfrak{sl}}
\newcommand{\fg}{\mathfrak{g}}
\newcommand{\ff}{\mathfrak{f}}
\newcommand{\fgl}{\mathfrak{gl}}
\newcommand{\fosp}{\mathfrak{osp}}
\newcommand{\fm}{\mathfrak{m}}

\newcommand{\ttau}{\tilde\tau}

\newcommand{\str}{\mathrm{str}}
\newcommand{\Sym}{\mathrm{Sym}}
\newcommand{\tr}{\mathrm{tr}}
\newcommand{\defi}{\mathrm{def}}
\newcommand{\Ber}{\mathrm{Ber}}
\newcommand{\spec}{\mathrm{Spec}}
\newcommand{\sschemes}{\mathrm{(sschemes)}}
\newcommand{\sschemeaff}{\mathrm{ {( {sschemes}_{\mathrm{aff}} )} }}
\newcommand{\rings}{\mathrm{(rings)}}
\newcommand{\Top}{\mathrm{Top}}
\newcommand{\sarf}{ \mathrm{ {( {salg}_{rf} )} }}
\newcommand{\arf}{\mathrm{ {( {alg}_{rf} )} }}
\newcommand{\odd}{\mathrm{odd}}
\newcommand{\alg}{\mathrm{(alg)}}
\newcommand{\sa}{\mathrm{(salg)}}
\newcommand{\sets}{\mathrm{(sets)}}
\newcommand{\SA}{\mathrm{(salg)}}
\newcommand{\salg}{\mathrm{(salg)}}
\newcommand{\varaff}{ \mathrm{ {( {var}_{\mathrm{aff}} )} } }
\newcommand{\svaraff}{\mathrm{ {( {svar}_{\mathrm{aff}} )}  }}
\newcommand{\ad}{\mathrm{ad}}
\newcommand{\Ad}{\mathrm{Ad}}
\newcommand{\pol}{\mathrm{Pol}}
\newcommand{\Lie}{\mathrm{Lie}}
\newcommand{\Proj}{\mathrm{Proj}}
\newcommand{\rGr}{\mathrm{Gr}}
\newcommand{\rFl}{\mathrm{Fl}}
\newcommand{\rPol}{\mathrm{Pol}}
\newcommand{\rdef}{\mathrm{def}}

\newcommand{\uspec}{\underline{\mathrm{Spec}}}
\newcommand{\uproj}{\mathrm{\underline{Proj}}}

\newcommand{\sym}{\cong}

\newcommand{\al}{\alpha}

\newcommand{\lam}{\lambda}
\newcommand{\de}{\delta}
\newcommand{\D}{\Delta}
\newcommand{\s}{\sigma}
\newcommand{\lra}{\longrightarrow}
\newcommand{\ga}{\gamma}
\newcommand{\ra}{\rightarrow}

\newcommand{\NOTE}{\bigskip\hrule\medskip}
%\copyrightnotice{<yyyy>}{<vol#>}
%{<start_page_num>}{<end_page_num>} %% year, volume, first page, last page.

%\setcounter{page}{<start_page_num>}

\title[Chiral Superspaces]
{The quantum chiral Minkowski and conformal superspaces}

%\arxurl{<hep_reference_#>}

\author[Cervantes, Fioresi, Lled\'{o}]
{D. Cervantes$^1$, R. Fioresi$^2$, M. A. Lled\'{o}$^3$}

\address{${}^1$Instituto de Ciencias Nucleares,
\\ Universidad Nacional
Aut\'{o}noma de M\'{e}xico \\
Circuito Exterior M\'{e}xico D.F. 04510, M\'{e}xico}
\address{
${}^2$Dipartimento di Matematica,  Universit\`{a} di
Bologna, \\ Piazza di Porta S. Donato, 5. 40126 Bologna. Italy.}
%lines should be separated with double backslashes: \\
\address{
${}^3$Departament de F\'{\i}sica Te\`{o}rica, \\
Universitat de Val\`{e}ncia and IFIC (CSIC-UVEG)
Fundaci\'{o} General Universitat de Val\`{e}ncia. \\ C/Dr.
Moliner, 50, E-46100 Burjassot (Val\`{e}ncia), Spain.}

\addressemail{daliac@nucleares.unam.mx, fioresi@dm.UniBo.it,
maria.lledo@ific.uv.es }

\begin{abstract}
We give a quantum deformation
of the  chiral super Minkowski space in four dimensions as the
big cell inside a quantum super Grassmannian. The quantization is performed in such way that the actions of the Poincar\'{e} and conformal quantum supergroups on the quantum Minkowski and quantum conformal superspaces are presented.
\end{abstract}

\maketitle

\section{Introduction}

The Minkowski space in four dimensions is just $\R^4$ with the
pseudoeuclidean metric $\rdiag(1,-1,-1,-1)$.
The Poincar\'{e} group $\rP(1,3)=\rSO(1,3)\propto \R^4$ is the group
that preserves such metric, while
the conformal group is the group that preserves the metric up to a global factor.
It is in fact the group $\rSO(2,4)$ and it acts non linearly on a
compactification of the Minkowski space, obtained by adjoining to it
not just a point at infinity, but the closure of a cone \cite{flv}. This compactification turns out to be the {\it Grassmannian
manifold} $G(2,4)$, that is, the space  of 2-planes inside a four dimensional
vector space and the Poincar\'{e} group together with the
dilations is precisely the subgroup of $\rSO(2,4)$ consisting of the
elements that leave the Minkowski space invariant.
The rest of the conformal transformations may send a point in the
Minkowski space to a point at infinity. The Minkowski space
sits inside the Grassmannian $G(2,4)$ as its
{\it big cell}, which is a dense open set inside it.
We will refer to the Grassmannian $G(2,4)$ as the {\it conformal space}.

The relation between the Poincar\'{e} and conformal group on one
side and the Minkowski and conformal space on the other side, are well
known, see for example Ref. \cite{pe}. A brief but complete review can
also be found in Ref. \cite{va}, which we have followed very closely
in spirit and notation. The analysis in there starts by considering the spin
group  of $\rSO(2,4)$, which is $\rSU(2,2)$, that contains the spin
group  of $\rSO(1,3)$, which is $\rSL(2,\C)_\R$. It is
more natural to work  in the complexified spaces  ($\rSL(4,\C)$ and $\rSL(2,\C)\times \rSL(2,\C)$ respectively) and to look at the end for the particular  real form associated to Minkowskian signature.

\medskip

This approach is very useful when extending the results to the
Minkowski and conformal {\it superspaces} (see \cite{flv}), since the
action on spinors is explicit in the formalism. The fascinating subject
of supergeometry emerges here as a very natural
framework. Supergeometry extends standard algebraic and differential
geometry in a less dramatic way than non commutative geometry in general \cite{co}. The
category of algebras considered in supergeometry are non commutative, but
their non commutativity affects only to some generators that {\it anticommute}. These are the  {\it odd generators}.

The {\it functor of points} is another tool, that one
borrows from standard algebraic geometry and extends to the
super setting. For a standard algebraic variety, the {\it geometric points} are the morphisms from the coordinate ring of the supervariety to the ground
field. If one considers the morphisms from the coordinate ring to another commutative ring, say $R$, one has the $R$-points of the algebraic variety. For a supervariety one takes $R$  to be a commutative superalgebra, and the $R$-points still make sense. The geometric points turn out to give the points of a standard variety called the {\it reduced variety}, and the odd variables  disappear. Only by considering morphisms to superalgebras instead than to commutative algebras or just the ground field one can recover the role of the odd variables, and with it all the information to reconstruct the supervariety.

We then see that via the functor of points the true nature of the odd variables appears
in a beautiful, perfectly consistent framework \cite{dm,va,cf}, which
is very close to the way in which physicists have been thinking and
talking about superspaces and supergroups. The terminology and perhaps
the level of rigor becomes more sophisticated, but a closer
look reveals many old concepts that have been used implicitly by physicists are at the core of the mathematical formulation.

\medskip

The next step is to produce a non commutative  version of the
Minkowski and conformal superspaces. In order to do this, we need to
substitute the commutative superalgebras by noncommutative ones, but
in that step the geometric intuition that we had retained in supergeometry with the functor of
points is lost. We have then to rely on the algebraic
counterpart of the geometric objects and try to generalize them to the non
commutative setting. Non commutative geometry
\cite{co} is certainly the most complete framework to do so, and
ultimately it will be connected to the quantization of space and superspace.

The approach that we follow here is the following: we substitute the supergroups by quantum supergroups and the corresponding homogeneous spaces by  quantum homogeneous spaces. This was the approach followed in Refs. \cite{fi2,fi3,fi4} for the non super case. We are then able to preserve also the realization of the quantized super Minkowski space as the `big cell' (appropriately defined) of the quantum conformal space.

 The problem of quantizing the Minkowski  superspace has appeared in many places in the physics literature. We mention some references, although our list is not exhaustive. We find a first step as early as in Ref.  \cite{sn} and \cite{bmn}, and deformations (mostly with constant Poisson bracket) inspired in string theory \cite{bgn,fel,se,flm}. In other papers one finds the quantization of a super Minkowski `phase space', which, although it is not exactly the problem that we examine here, it is also of interest \cite{bdpz,hrs,ykn}. Quantizations for other superspaces such as supergroups \cite{ma2,fi5,fi6}, their coadjoint orbits \cite{fl2} or other homogeneous superspaces $\C\mathrm{P}^{m|n}$ \cite{imt}  are constructed.

 The principle that guides us in choosing a particular deformation is
that we want them to preserve the action of the corresponding symmetry
groups and we also ask that the quantum Minkowski superspace appears
as the big cell (appropriately defined in algebraic terms) inside the quantum conformal
superspace. Obviously, all of these requests have to  be made precise in
the framework of deformation quantization.

\medskip

The content of the present paper is as follows.

We start in Section \ref{chiralsuperfields} by giving an overview
on the physical question that originates our discussion. We
are after a rigorous  mathematical
description of the chiral superfields and their quantization,
making sure to preserve the natural supergroup actions of
the superconformal and super Poincar\'{e} groups.

We will devote Section \ref{secfunctorofpoints} to give an intuitive
explanation of how the functor of points works in terms of superalgebras.

In Section \ref{secplucker} we review briefly the classical theory
of the Grassmannian manifold as embedded in the projective space.
This is the standard {\it Pl\"{u}cker embedding}. We then generalize these structures to the super setting.

In Section \ref{secquantumsupergrass} we
use the technology developed in the series of
papers \cite{fi2,fi3,fi4}, to quantize the Grassmannian and flag
supervarieties by replacing the symmetry supergroup by a quantum supergroup.
We then discuss quantum supergroups and their homogeneous
spaces \cite{ma2,fi5,fi6,flv}, by looking at the corresponding non
commutative superalgebras.
One result of this approach is that although the algebras become non
commutative,
the group law, represented by the {\it comultiplication} in the quantum
supergroup, is not deformed. This can be interpreted
by saying that the physical symmetry principle remains
intact in the process of quantization.

In Section \ref{secstarproduct} we give a definition of the {\it quantum big cell} inside the quantum Grassmannian, presenting a coaction of the quantum super Poincar\'{e} group.

In Section \ref{secconclusions} we state our conclusions
and give the guidelines for further work.

We have left for Appendix
\ref{supergeometry-app} a brief and formal account of
the most fundamental concepts of supergeometry.

%\paragraph
{\bf Notation.}

We will say that a superalgebra is {\sl commutative} if for two
elements
$a$ and $b$ of definite parities $p_a$ and $p_b$ we have
$$a\cdot b=(-1)^{p_ap_b}b\cdot a.$$
In physics this is usually called an {\sl (anti)commutative} or supercommutative algebra,
but we prefer to keep the word commutative as in \cite{dm}.

In this paper the word `classical' is used for commutative
superalgebras,
as opposed to the `quantum' superalgebras, which are noncommutative.

Also, we stress that all our definitions are done in terms of algebras
and superalgebras, since it is the only aspect that survives in the noncommutative case.
$\hbs$

\section{Real and chiral superfields in Minkowski superspace}
\label{chiralsuperfields}

We want to devote this section to introduce real and chiral superfields as they are used in physics as well as to motivate the importance of having them quantized. For this purpose we will consider superfields in super Minkowski space: one can then introduce the notion of conformal or superconformal invariance in quuantum field theory.

\subsection{Definitions}

We consider the complexified Minkowski space $\C^4$. The $N=1$ {\sl scalar superfields} on the complexified Minkowski space are elements of the commutative superalgebra  \be\cO(\C^{4|4})\equiv C^\infty(\C^4)\otimes \Lambda[\theta^1,\theta^2,\bar\theta^1,\bar\theta^2],\label{superminkowski}\ee where $ \Lambda[\theta^1,\theta^2,\bar\theta^1,\bar\theta^2]$ is the Grassmann (or exterior) algebra generated by the odd variables $\theta^1,\theta^2,\bar\theta^1,\bar\theta^2$.
Giving this superalgebra is equivalent to giving the superspace $\C^{4|4}$ as defined in  Appendix \ref{apbasicdefinitions}.

We will denote the coordinates (or generators) of the superspace
as
\begin{eqnarray*}&x^\mu, \qquad &\mu=0,1,2,3 \quad \hbox{(even
coordinates),}\\&\theta^\alpha,\bar \theta^{\dot\alpha},\qquad&
\alpha, \dot \alpha=1,2 \quad \hbox{(odd coordinates),}
\end{eqnarray*}
and a  superfield, in terms of its {\sl field components}, as
\begin{align*}\Psi(x, \theta,\bar \theta)=&\psi_0(x)+\psi_\alpha(x)\theta^\alpha+\psi'_{\dot\alpha}(x)\bar\theta^{\dot\alpha}+
\psi_{\alpha \beta} (x)\theta^{ \alpha} \theta^{ \beta}+\psi_{ \alpha\dot\beta} (x) \theta^{\alpha}\bar\theta^{\dot\beta}+\\&
\psi'_{\dot\alpha\dot\beta}(x)\bar\theta^{\dot\alpha}\bar\theta^{\dot\beta}+
\psi_{ \alpha \beta\dot\gamma}(x) \theta^{ \alpha} \theta^{ \beta}\bar\theta^{\dot\gamma}+
\psi'_{ \alpha\dot\beta\dot\gamma}(x) \theta^{ \alpha}\bar\theta^{\dot\beta}\bar\theta^{\dot\gamma}+
\psi_{ \alpha \beta\dot\gamma\dot\delta}(x) \theta^{ \alpha} \theta^{ \beta}\bar\theta^{\dot\gamma}
\bar\theta^{\dot\delta}.\end{align*}

A {\it conjugation} on a superalgebra $\cA$ (not necessarily commutative) is an antilinear involutive map satisfying
\be(f\cdot g)^*= (-1)^{p_fp_g}g^*f^*,\qquad f,g\in \cA\label{complexconj}\ee where $p_f$ is the parity of the element $f$. Here we take the convention of Ref. \cite{dm} (see page 89 in there for a detailed explanation) because it has categorical meaning. Moreover, it is much more appropriate for superalgebras that are not necessarily commutative. It differs from the one used in physics, because if  $f$ and $g$ are odd elements, then
$$(f\cdot g)^*=-g^*f^*,$$ and if the superalgebra is commutative
$$(f\cdot g)^*=f^*g^*.$$ Instead, physicists use an operation such that
$$(f\cdot g)^\rho= g^\rho f^\rho,$$ but this gives rise to minus signs and `\,i\,' factors that we would like to avoid.
As explained in  Ref. \cite{dm}, one convention can be reverted into the other by a change
$$f^*=\left\{\begin{array}{rl}f^\rho& \hbox {if $f$ is even,}\\\\\ri f^\rho& \hbox {if $f$ is odd,}\end{array}\right.$$ so
$$(f\cdot g)^\rho=g^\rho f^\rho.$$

On the space of complex functions  $C^\infty(\C^4)$, there exists the standard complex conjugation, denoted as   $$f^*=\bar f,\qquad f\in C^\infty(\C^4). $$ Hence, to give a conjugation it is enough to give it on the odd  generators. This is done formally in the following way
\be(\theta^\alpha)^*=
\bar\theta^{\dot\alpha},\qquad (\bar\theta^{\dot \alpha})^*=
\theta^{\alpha}.\label{cc}\ee This is then extended by (\ref{complexconj}) and antilinearity to the whole superalgebra $\cO(\C^{4|4})$. Real superfields then belong to $\cO(\R^{4|4})$.

Using the convention (\ref{complexconj}) a superfield  $\Psi(x,\theta,\bar \theta)$ is {\it real} if and only if its field components satisfy
\begin{align*}&\psi_0^*=\psi_0,\qquad \psi_\alpha^*=\psi'_{\dot\alpha},\qquad \psi_{\alpha\beta}^*=\psi'_{\dot\alpha\dot\beta} \qquad\psi_{\alpha\dot\beta}^*=-\psi_{\beta\dot\alpha},\\ &\psi_{\alpha\dot\beta\dot\gamma}^*=\psi'_{\gamma\dot\alpha\dot\beta},\qquad
\psi_{\alpha\beta\dot\gamma\dot\delta}^*=\psi_{\gamma\delta\dot\alpha\dot\beta}.\end{align*}

%\paragraph
{\bf Action of the Lorentz group SO(1,3).} There is an action of the double covering of the complexified Lorentz group,
$\rSpin(1,3)^c\approx\rSL(2,\C)\times \rSL(2,\C)$ over $\C^{4|4}$. The even coordinates  $x^\mu$ transform
in the fundamental representation of $\rSO(1,3)$ ($V$),
$$x^\mu\mapsto \Lambda^\mu{}_\nu x^\nu,$$
while $\theta$ and $\bar \theta$ are Weyl spinors (or
half spinors). More precisely, the coordinates $\theta$ transform in one of the spinor representations, say
$S^+\approx(1/2,0)$, and $\bar \theta$ transform in the opposite chirality representation, $S^-\approx(0,1/2)$,
$$\theta^\alpha\mapsto {S}^\alpha{}_\beta\theta^\beta,\qquad \bar\theta^{\dot\alpha}\mapsto \tilde S^{\dot \alpha}{}_{\dot \beta}\theta^{\dot \beta}.$$
In fact, for the real form $\rSpin(1,3)$, the representations
$S^+$ and $S^-$ are complex, and they are related by complex
conjugation, so this is consistent with the rule (\ref{cc}).

 The scalar superfields are invariant under the action of the Lorentz group,
$$\Psi(x,\theta,\bar\theta)= (R\Psi)(\Lambda^{-1} x, S^{-1}\theta,\tilde S^{-1}\bar\theta),$$ where $R\Psi$ is the superfield obtained by transforming the  the field components
$$R\psi_0(x)=\psi_0(x),\qquad R\psi_\alpha(x)=S_\alpha{}^\beta\psi_\beta(x), \;\; \dots$$

The hermitian matrices
$$\sigma^0=\begin{pmatrix}1&0\\0&1\end{pmatrix},\quad
\sigma^1=\begin{pmatrix}0&1\\1&0\end{pmatrix},\quad\sigma^2=\begin{pmatrix}0&-\ri\\\ri&0\end{pmatrix},\quad
 \sigma^3=\begin{pmatrix}1&0\\0&1\end{pmatrix},$$ define a
  $\rSpin(1,3)$-morphism
$$\begin{CD}S^+\otimes S^-@>>>V\\
s^\alpha \otimes  t^{\dot\alpha} @>>>
s^\alpha\sigma^\mu_{\alpha\dot\alpha}t^{\dot\alpha}.
\end{CD}$$$\hbs$

% \paragraph
{\bf Derivations.}
A {\sl left derivation} of degree $m=0,1$ of  a super algebra $\cA$ is a
linear map $D^L:\cA\mapsto \cA$ such that
$$
D^L(\Psi\cdot\Phi)=D^L(\Psi)\cdot \Phi +(-1)^{mp_\Psi}\Psi\cdot D^L(\Phi).
$$
Graded left derivations span a $\Z_2$-graded vector space (or {\sl supervector space}).

In general, linear maps over a supervector space are also  a $\Z_2$-graded vector space. A map has degree 0 if it preserves the
parity and degree 1 if it changes the parity. For the case of derivations of a commutative superalgebra,
 an even derivation has degree 0 as a linear map and an odd derivation
has degree 1 as a linear map.

In the same way one defines  {\sl right derivations},
$$
D^R(\Psi\cdot \Phi)=(-1)^{mp_\Phi}D^R(\Psi)\cdot \Phi +\Psi\cdot D^R(\Phi).
$$
Notice that  derivations of degree zero are both, right and left derivations.
Moreover, given a left derivation $D^L$ of degree $m$  one can define a
right derivation $D^R$ also  of degree $m$
 in the following way
\be D^R \Psi=
(-1)^{m(p_\Psi+1)}D^L\Psi.\label{rl}\ee

Let us now focus on the commutative superalgebra $\cO(\C^{4|4})$. We define the standard left derivations
\begin{align*}\partial_\alpha^L\Psi&=\psi_\alpha+
2\psi_{\alpha \beta} \theta^{ \beta}+\psi_{ \alpha\dot\beta} \bar\theta^{\dot\beta}+
2\psi_{ \alpha \beta\dot\gamma}  \theta^{ \beta}\bar\theta^{\dot\gamma}+
\psi'_{ \alpha\dot\beta\dot\gamma} \bar\theta^{\dot\beta}\bar\theta^{\dot\gamma}+
2\psi_{ \alpha \beta\dot\gamma\dot\delta}  \theta^{ \beta}\bar\theta^{\dot\gamma}
\bar\theta^{\dot\delta},\\
\partial^L_{\dot\alpha}\Psi&=\psi'_{\dot\alpha}-\psi_{ \beta\dot\alpha} \theta^{\beta}+
2\psi'_{\dot\alpha\dot\beta}\bar\theta^{\dot\beta}+
\psi_{ \gamma \beta\dot\alpha} \theta^{ \gamma} \theta^{ \beta}\bar-2
\psi'_{ \beta\dot\alpha\dot\gamma} \theta^{ \beta}\bar\theta^{\dot\gamma}+2
\psi_{ \gamma \beta\dot\alpha\dot\delta} \theta^{ \gamma} \theta^{ \beta}\bar
\bar\theta^{\dot\delta}.
\end{align*}
With our convention (\ref{complexconj}), one has that
$$(\partial^L_\alpha\Psi)^*=\partial^L_{\dot\alpha}\Psi^*.$$
Also, using (\ref{rl})  one can define $\partial_\alpha^R, \partial^R_{\dot\alpha}$. They have the same property than the left derivatives under complex conjugation.

We consider now the odd left derivations
$$Q^L_\alpha=\partial^L_\alpha
-\ri\sigma^\mu_{\alpha\dot\alpha}\bar\theta^{\dot\alpha}\partial
_\mu,\qquad  \bar Q^L_{\dot \alpha}=-\partial^L_{
\dot \alpha}
+\ri\theta^{\alpha}\sigma^\mu_{\alpha\dot\alpha}\partial_\mu.$$ They satisfy the anticommutation rules
$$\{Q^L_\alpha, \bar Q^L_{\dot \alpha}\}=2\ri \sigma^\mu_{\alpha\dot\alpha}\frac\partial{\partial
x^\mu}, \qquad \{Q^L_\alpha, Q^L_{\beta}\}=\{\bar Q^L_{\dot
\alpha}, \bar Q^L_{\dot \beta}\}=0.$$ $Q^L$ and $\bar Q^L$ are the
{\sl supersymmetry charges} or {\sl supercharges}. Together with
$$P^\mu=-\ri \frac\partial{\partial
x^\mu},$$ they form a Lie superalgebra,  the {\sl supertranslation
algebra}, which then acts on  the superspace $\C^{4|4}$.

Let us define another set of (left) derivations,
$$D^L_\alpha=\frac{\partial^L}{\partial \theta^\alpha}
+\ri\sigma^\mu_{\alpha\dot\alpha}\bar\theta^{\dot\alpha}\frac\partial{\partial
x^\mu},\qquad  \bar D^L_{\dot \alpha}=-\frac{\partial^L}{\partial
\bar\theta^{\dot \alpha}}
-\ri\theta^{\alpha}\sigma^\mu_{\alpha\dot\alpha}\frac\partial{\partial
x^\mu},$$ with anticommutation rules
$$\{D^L_\alpha, \bar D^L_{\dot \alpha}\}=-2\ri \sigma^\mu_{\alpha\dot\alpha}\frac\partial{\partial
x^\mu}, \qquad \{D^L_\alpha, D^L_{\beta}\}=\{\bar D^L_{\dot
\alpha}, \bar D^L_{\dot \beta}\}=0.$$ They also form a Lie
superalgebra, isomorphic to the supertranslation algebra. This can be seen by taking
$$Q^L\rightarrow -D^L,\qquad \bar Q^L\longrightarrow \bar D^L.$$

It is easy to see that the supercharges anticommute with the
derivations $D^L$ and $\bar D^L$. For this reason, $D^L$ and $\bar
D^L$ are called {\sl supersymmetric covariant derivatives} or
simply {\sl covariant derivatives}, although they are not related
to any connection form.

$\hbs$

We go now to the central definition.

\begin{definition}

A {\sl chiral superfield} is a superfield $\Phi$ such that
\begin{equation}\bar D^L_{\dot\alpha}\Phi=0.\label{chiralcons}\end{equation}$\hbs$\end{definition} Because of the anticommuting
properties of $D's$ and $Q's$, we have that
$$ \bar D^L_{\dot\alpha}\Phi=0\quad \Rightarrow\quad  \bar
D^L_{\dot\alpha}(Q^L_\beta\Phi)=0,\quad \bar D^L_{\dot\alpha}(\bar
Q^L_{\dot\beta}) \Phi=0.$$ This means that the supertranslation
algebra acts on the space of chiral superfields.

On the other hand, due to the derivation property,
$$\bar D^L_{\dot\alpha}(\Phi\Psi)=\bar D^L_{\dot\alpha}(\Phi)\Psi+(-1)^{p_\Phi}\Phi\bar
D^L_{\dot\alpha}(\Psi),$$ we have that the product of two chiral
superfields is again a chiral superfield.

\subsection{Shifted coordinates}
One can solve the constraint (\ref{chiralcons}) in the following way. Notice that the
quantities
\be y^\mu=x^\mu
+\ri\theta^\alpha\sigma^\mu_{\alpha\dot\alpha}\theta^{\dot\alpha},\qquad
\theta^\alpha\label{shifted}\ee
satisfy
$$\bar D^L_{\dot\alpha}y^\mu=0,\qquad \bar
D^L_{\dot\alpha}\theta^\alpha=0,$$ and using the derivation
property, any superfield of the form
$$\Phi(y^\mu,\theta),\qquad \hbox{satisfies  }\quad \bar D^L_{\dot\alpha}\Phi=0$$ and so it  is a chiral
superfield. This is the general solution of (\ref{chiralcons}).

We can make the change of coordinates
$$x^\mu, \; \theta^\alpha,\; \bar\theta^{\dot\alpha}\;\longrightarrow \; y^\mu=x^\mu +i\theta^\alpha
\sigma^\mu_{\alpha\dot\alpha} \bar\theta^{\dot \alpha},\;
\theta^\alpha,\;\bar\theta^{\dot\alpha}.$$

A superfield may be expressed in both coordinate systems
$$\Phi(x, \theta, \bar\theta)=\Phi'(y, \theta, \bar\theta).$$ The
covariant derivatives and supersymmetry charges take the form
\begin{eqnarray*}
D^L_\alpha\Phi'  =\frac{\partial^L \Phi'}{\partial\theta^\alpha}
 +2i
\sigma^\mu_{\alpha\dot\alpha}\bar\theta^{\dot\alpha}
\frac{\partial^L \Phi' }{\partial y^\mu}\qquad \bar
D^L_{\dot\alpha}\Phi' =
-\frac{\partial^L \Phi'}{\partial\bar\theta^{\dot\alpha}},\\
  \bar Q^L_{\dot\alpha}\Phi'  =
-\frac{\partial^L \Phi'}{\partial\bar\theta^{ \dot\alpha}} +2i
\theta^\alpha\sigma^\mu_{\alpha\dot\alpha}\frac{\partial^L \Phi'
}{\partial y^\mu}\qquad  Q^L_\alpha\Phi' =
\frac{\partial^L\Phi'}{\partial\theta^\alpha}.
\end{eqnarray*}
In the new coordinate system the chirality condition is simply
$$\frac{\partial^L\Phi'}{\partial\bar\theta^{\dot\alpha}}=0,$$ so it
is similar to a holomorphicity condition on the $\theta$'s.

This shows that chiral scalar superfields are elements of the commutative superalgebra $\cO(\C^{4|2})=\C^\infty(\C^4)\otimes \Lambda[\theta^1,\theta^2]$.
We shall realize this superspace as the big cell inside the chiral conformal superspace, which is the Grassmannian of $2|0$-subspaces of $\C^{4|1}$.

The complete (non chiral) conformal superspace is in fact the flag supervariety of $2|0$-subspaces inside $2|1$-subspaces of  $\C^{4|1}$. On this supervariety one can put a reality condition, and the real Minkowski superspace is the big cell inside the superflag.
It is instructive to compare Eq. (\ref{shifted}) with the incidence relation for the big cell of the flag manifold in Eq. (12) of Ref. \cite{flv}. We can then be convinced that the Grassmannian that we use to describe chiral superfields is inside the (complex) superflag.

\subsection{Supersymmetric theories}
Wess-Zumino models are supersymmetric models for one or several chiral superfields. These were the first type of supersymmetric theories that were written down \cite{wz1}. Chiral superfields also appear in super Yang-Mills theories \cite{wz2,fz}, where the parameter of the gauge transformation is itself a chiral superfield.

The study of the chiral super Minkowski space is then justified from the physical point of view. Of course, most of the theories make use of the real super Minkowski space, and one needs also to consider real fields to formulate supersymmetric theories.

Also, it is important to consider the embedding of super Minkowski space inside conformal superspace, since some theories (for example, some Wess-Zumino models and N=4 super Yang-Mills theory) have this symmetry.
In fact, for the classic (non quantum) case, this has been done in the modern language of supergeometry in \cite{flv}.

In this paper though, we want to consider a quantization of these superspaces that preserves the action of the corresponding supergroups. It has been particularly difficult to find deformations of the space of chiral superfields involving also the odd variables \cite{fel,flm}. Up to now, this has prevented to formulate Wess-Zumino or Yang-Mills models in a non commutative superspace with a non trivial deformation of the odd part and preserving the supersymmetry. Essentially, what happened in previous formulations is that the covariant derivatives were not anymore derivations of the noncommutative product, and then the ring of chiral superfields did not extend to a quantum chiral ring. Some proposals to keep a chiral ring (but not an antichiral one) include the partial (explicit) breaking of supersymmetry \cite{se,fel}.

In our formulation, we start with the classical chiral ring and find a {\it quantum chiral ring} in a natural way. We substitute the supergroup by a quantum supergroup and preserving the relations among all the elements of the construction. As it is well known, the comultiplication is not deformed when going from the classic to the quantum group, which means that the supersymmetry algebra is preserved without deformation although now it is realized on a non commutative superspace. Mathematically, this is already a non trivial problem, and physically it is a problem that must be solved in order to formulate the theories that use chiral superfields in non commutative spaces.

Our approach will be complete once we extend it to the real super Minkowski space. In order to do this, one has to deal with the flag supermanifold. The Grassmannian then sits inside the complexified flag supermanifold.

 One can certainly extend the same philosophy of quantization to the flag supervariety. Nevertheless, the problem is non trivial, presents its own complications and will be the subject of a forthcoming paper.

Finally, since the superconformal symmetry is implicit in our approach, we expect to obtain in the future a basis to formulate conformal theories in a non commutative space. This will include for example $N=4$ super Yang-Mills.

\section{The functor of points}\label{secfunctorofpoints}

The functor of points is an extremely useful tool in classical
algebraic geometry, which becomes essential in supergeometry in order
to recover the geometric intuition, otherwise lost.
In this section we shall give an intuitive, though operative, summary
of its definitions and properties, sending
the reader to the Appendix \ref{supergeometry-app} for the complete
treatment and all the references.

 For definiteness, we take the ground field to be  $k=\R,\C$. A {\it
superalgebra} $\cA$ is a $\Z_2$-graded algebra, $\cA=\cA_0 \oplus \cA_1$,
where  $p(x)$ denotes the parity of a homogeneous element $x$,
so $p(x)=0$ if $x\in \cA_0$ and $p(x)=1$ if $x\in \cA_1$. The subspace
$\cA_0$ is an algebra, while the subspace $\cA_1$ is an $\cA_0$-module.

The superalgebra $\cA$ is said to be {\it commutative} if for any two
homogeneous elements $x, y$
$$
xy=(-1)^{p(x)p(y)}yx.
$$From now on all superalgebras are assumed to be commutative
unless otherwise specified. The category of commutative superalgebras will be denoted by $\salg$, and the category of sets is denoted by (sets).

A {\it functor}
$h:\mathrm{(salg)}\longrightarrow \mathrm{(sets)}$
is {\it representable} if there exists a superalgebra $\cB$ such that
\be
\begin{CD}
h:\mathrm{(salg)} @>>>\sets \\
\cA @>>>h(\cA)=\Hom(\cB, \cA).
\end{CD}\label{affinefunctor}\ee
 We say that the superalgebra representing the functor $h$ is $\cB$. In that case, we will denote the functor $h$ as $h_\cB$, to stress the fact that it is represented by the superalgebra $\cB$. An element of $\Hom(\cB,\cA)$ is called an {\it  $\cA$-point of $h_\cB$}.

If we restrict the functor $h$ to the category of $k$-algebras
(that is, if we remove the odd generators), and we also demand that the algebras are {\it reduced}\footnote{An algebra is reduced if it has no nilpotent elements.} and finitely generated (these algebras are often called {\it affine algebras}), then a representable functor
$h$ corresponds to an algebra $\cO(X)$,
that is  the coordinate ring of an {\it affine variety} $X$.

\begin{example}
Let $S^2$ be the unit 2-sphere in $\C^3$. Its coordinate ring
is  given by
$$
\cO(S^2)=\C[x,y,z]/(x^2+y^2+z^2-1).
$$
The elements in $\cO(S^2)$ correspond to the polynomial
functions on the variety $S^2 \subset \C^3$.
Each morphism $\cO(S^2)\rightarrow\C$ is given in terms of the
images of the generators
$$
(x,y,z)\longrightarrow (a,b,c)\; \hbox{  with  } a,b,c\in  \C
\hbox{  such that  } a^2+b^2+c^2=1,
$$
hence, any such morphism represents a point of the sphere.
If we consider morphisms $\cO(S^2)\rightarrow \mathbf Q$, they are in
one-to-one correspondence with rational points on the sphere.
In this example the geometric points (that is the $\C$-points)
allow us to recover all the information on $S^2$.
$\hbs$
\end{example}

In Appendix \ref{supergeometry-app} we give the definition of
affine supervariety, however, for all the practical purposes in
this paper, we can identify an affine supervariety $X$ with its
coordinate ring $\cO(X)$, this time $\cO(X)$ being a superalgebra. As in the non super case, it results that affine supervarieties are in one to one correspondence with  {\it affine superalgebras}. This means only that  the ordinary  algebra
obtained by taking modulo by the odd ideal $\cO(X)/\cI_{\mathrm{odd}}$
is an affine algebra. Then, there is an affine variety corresponding to
$\cO(X)/\cI_{\mathrm{odd}}$. This is an affine variety underlying
the affine supervariety.

The $\cA$-points of the supervariety $X$ are  the
morphisms $\cO(X) \lra \cA$ and all
the information about the affine supervariety is fully encoded
in its coordinate superalgebra, or equivalently in its functor of
points $$\begin{CD}h_X:\salg @>>>  \sets \\\cA @>>>h_X(\cA)=\Hom(\cO(X),\cA).\end{CD}$$
We have denoted $h_X$ instead of $h_{\cO(X)}$ for simplicity,
since there is no possibility of confusion in this context.

To clarify these concepts we describe a simple example.

\begin{example} \label{affinefopts}
 We define the
\textit{polynomial superalgebra} as:
$$
k[x^1,\dots ,x^p,\theta^1,\dots ,\theta^q]
:= k[x^1,\dots ,x^p]\otimes \Lambda[\theta^1,\dots ,\theta^q]
$$
We want to interpret this superalgebra as
 the coordinate ring of the affine superspace of
superdimension $p|q$,
that we shall denote with the symbol $k^{p|q}$.
If $\cA$ is a generic (commutative) superalgebra, an $\cA$-point of
$k^{p|q}$ is given by a morphism
$k[x^1,\dots ,x^p,\theta^1,\dots ,\theta^q] \lra \cA$, which is
determined once we know the images of the generators $$(x^1,\dots ,x^p,\theta^1,\dots ,\theta^q)\lra(a^1 \dots a^p,\al^1 \dots \al^q),$$
with $a^i \in \cA_0$ and $\al^j \in \cA_1$. Notice that the $k$-points
of $k^{p|q}$ are given by $(k_1 \dots k_p, 0 \dots 0)$ and coincide
with the points of the affine space $k^p$. From this example it is
clear that the knowledge of the geometric points, that
is the $k$-points, is by no means
sufficient to describe the supergeometric object.

$\hbs$
\end{example}

There is an important property of the functor of points, which
is constantly used and makes all calculations  easier.
In complete analogy with the classical setting,
the functor of points of an affine supervariety
is determined by its image on the {\it local
superalgebras}\footnote{A local superalgebra is a superalgebra  that has a unique,
maximal ideal which is homogeneous with respect to the
$\Z_2$ grading.} (For the proof of this fact see Ref. \cite{cf} ch. 10).
Local superalgebras are important, since their behavior resembles
those of fields, thus rendering the functor of points easier to
describe in the examples that we are interested in.

\medskip

We now turn to projective supervarieties.
In the classical setting projective varieties are harder to describe,
since the homogeneous coordinate ring that we associate to a projective
variety encodes not only the structure of the variety, but also
its embedding into projective space. Such ring hence
has more information than the variety itself. As a consequence we have that non isomorphic coordinate homogeneous
rings may correspond to the same projective variety, a phenomenon
that we do not see in the affine case, where coordinate rings and
affine varieties correspond bijectively to each other and contain essentially
the same information. A fancy way to express this, is to say that
the category of affine (super)varieties over $k$ is equivalent to the category of affine
$k$-(super)algebras. So in the affine setting we have that in both, the
classical and the super setting,
$$
\Hom(\cO(X),\cO(Y))=\Hom(Y,X).
$$
Notice that the role of $X$ and $Y$ are interchanged when passing
from the coordinate (super)rings to the (super)varieties.

There are two equivalent, but different, ways to approach projective
varieties in the ordinary setting and we shall briefly describe them.
Both can be generalized to the super setting. For us the second one is far more important, since it gives a setting suitable for the quantization. The first one will also be used (implicitly) when we perform the quantization of the big cell in Sec. \ref{secstarproduct}.

One way to approach
projective varieties is to view them locally as affine varieties
that can suitably  be patched together. In other words a projective
variety is a topological space that is
covered by affine varieties, whose coordinates, in the overlaps
of different affine varieties,
behave in a certain way.
We can build the functor of points of a projective variety by giving
the functors of points of the affine varieties and then asking that they satisfy certain gluing conditions. This is in essence
the meaning of the Representability Theorem\footnote{The
Representability Theorem is stated in the more
general setting of superschemes, which are a far reaching generalization
of projective supervarieties.}
in Appendix
\ref{supergeometry-app}.

A second equivalent way to define a projective variety $X$ is
to look at the points in a projective space $\bP^n$, satisfying
homogeneous equations $f_1=0,
\dots f_n,=0$ (homogeneous here refers to the $\Z$-grading of polynomials).
Hence the homogeneous ring $S=k[x_0 \dots  x_n]/(f_1 ,\dots, f_n)$
determines uniquely the projective variety $X$ together with its
embedding into $\bP^n$. However, as we have already remarked, the
variety $X$ does not determine uniquely the
ring $S$:
there could in fact
be a different embedding of $X$ into some other projective space, yielding
a homogeneous algebra non isomorphic to $S$. For this reason the
functor of points of a projective variety is more tricky and
it is not directly related to the coordinate ring $S$ of the variety
itself as it happens for the affine case. In particular this
functor is {\sl never} representable, hence making the theory more
difficult since there is no algebra that can be naturally associated to it.
Hence we are facing a new problem: how can we decide whether
a functor $h:\salg \lra \sets$ is the functor of points of a
projective supervariety? Certainly we cannot say as before that
this is equivalent for $h$ to be representable. The answer is again
the Representability Theorem in Appendix \ref{supergeometry-app}.

As an example of its application,
let us examine the functor of points of
the projective space $\bP^n$.

\begin{example}
Let us consider the functor:
$h: \salg \lra \sets$, where $h(\cA)$ are the projective $\cA$-modules of
rank one in $\cA^n$

Equivalently $h(\cA)$ consists
of the pairs $ (L,\phi)$, where $L$ is a projective
$\cA$-module of rank one, and $\phi$ is a surjective morphisms $\phi: \cA^{n+1} \lra L$. These pairs are taken modulo the equivalence relation
$$(L,\phi)\approx(L',\phi')\qquad \Leftrightarrow\qquad L\approx L',\qquad \phi'=a\circ\phi,$$
where $a:L\rightarrow L'$ is the isomorphism. We can substitute the last condition by saying that $\phi$ and $\phi'$ have the same kernel.

 If $\cA=k$ is a field, then projective modules are free and a morphism
$$\phi :k^{n+1}\rightarrow k$$
is specified by a $n$-uple $(a^1,\dots a^{n+1})$, with $a^i\in k$, not all of the $a^i=0$. The equivalence relation becomes
$$(a^1,\dots a^{n+1})\sim(b^1,\dots b^{n+1}) \quad \Leftrightarrow\quad (a^1,\dots a^{n+1})=\lambda(b^1,\dots b^{n+1}),$$
with $\lambda\in k^\times$ understood as an automorphism of $k$. It is clear then that
$h(k)$ consists of all the lines through the origin
in the vector space $k^{n+1}$.

If $\cA$ is local, projective modules are free over local
rings. We then have a situation similar to the field setting:
equivalence classes are lines in the $\cA$-module $\cA^{n+1}$.
The functor $h$ is `covered' by  functors which on a local algebra $\cA$
$$h_i(\cA)=\left\{\rspan\{(a^1, \dots, a^{n+1})\}\in\cA^{n+1}\;|\; a^i \hbox{ is invertible}\right\}.$$
In fact, in order to correspond
to a surjective morphism, the $n$-uple $(a^1, \dots ,a^{n+1})$ must contain
at least an invertible $a^i$. Hence, using the
equivalence relation, we have that
$h_i(\cA)=\{(a_1/a_i, \dots ,\hat a_i,\dots a^{n+1}/a_i)\}$ and consequently $h_i=h_{\bA^n}$,
as described in the Example \ref{affinefopts}.

Using the Representability Theorem one can then show that the functor $h$ is the functor
of points of a variety that we call the projective space and whose
geometric points  coincide with the projective
space $\bP^n$ over the field $k$ as we usually understand it.
$\hbs$
\end{example}

This example can be easily generalized to the supercontext:
we consider the functor $h_{\bP^{m|n}}:\salg \lra \sets$, where
$h_{\bP^{m|n}}(\cA)$ is defined as the set
the projective $\cA$-modules of rank $1|0$ in $\cA^{m|n}:=\cA \otimes k^{m|n}$.
This is  the functor of points of a supervariety,
called the {\it projective superspace}.

\medskip

The next question that  we want to tackle is how we can define an
embedding of a (super)variety into the projective (super)space
using the functor of points notation.

Let $X$ be a projective supervariety and
$\Phi:X \lra \bP^{n|m}$ be an injective morphism. In the notation of
the functor of points, $\Phi$ is a {\it natural transformation} between the two functors $h_X$ and $h_{\bP^{n|m}}$, given by  $$\Phi_{\cA}: h_X(\cA) \lra h_{\bP^{n|m}}(\cA)$$
with $\Phi_{\cA}$ injective.

If $\cA$ is a local
superalgebra, then an $\cA$-point in $h_{\bP^{n|m}(\cA)}$
is in $\phi_\cA(h_X(\cA))$
if and only if it satisfies certain homogeneous polynomial relations
$f_1=0, \dots ,f_r=0$ in the indeterminates $x_1, \dots x_n,\xi_1, \dots \xi_m$.
Moreover $X$ is the projective
supervariety associated to
$k[x_1, \dots, x_n,\xi_1, \dots, \xi_m]/(f_1, \dots, f_r)$.
In other words the superalgebra
$k[x_1 \dots x_n,\xi_1 \dots \xi_m]/(f_1 \dots f_r)$ is the coordinate
superalgebra of the projective super variety $X$ with respect to
the embedding $\phi$. Recall  that,
although $X$ does not determine a
coordinate superalgebra as in the affine case,
 the coordinate superalgebra of $X$
with respect to an embedding into projective superspace
does determine the projective supervariety $X$.

In summary, to determine the coordinate superalgebra of a
projective supervariety with
respect to a certain projective embedding, we need
to check the relations satisfied by the coordinates {\sl just on local
superalgebras}. This will be our starting point in Section
\ref{pluckersupergrass}
when we shall determine the coordinate superalgebra of the Grassmannian
supervariety with respect to its Pl\"{u}cker embedding.

\medskip

As a side remark we want to notice that the Grassmannian
supervariety is not in general a projective supervariety \cite{kn,ma1}, contrary to the
classical setting. This unpleasant feature luckily does not affect
us, since the particular Grassmannian variety
that we see in
this paper is projective.

\section{The super Grassmannian variety}
\label{secplucker}
We will start by reviewing the Pl\"{u}cker embedding of
the Grassmannian variety. We then turn to a description of the
same object as a {\sl quotient space} for the action of the
special linear supergroup and we give an
explicit description of the big cell inside the Grassmannian
supervariety. The big cell is especially important since it is
identified with the complex Minkowski superspace, while the subgroup
of $\rSL(4|1)$ stabilizing the big cell contains
the Poincar\'{e} supergroup times  the dilations.

Our point of view is purely algebraic, since it is
the most suitable for the quantization. The reader however must be aware that
these objects have a natural differential structure, as always is
the case for smooth supergroups and the supervarieties on which
they are acting \cite{fi7}.

\subsection{Pl\"{u}cker embedding of the Grassmannian variety
and the big cell\label{plucker}}

The Grassmannian variety $G(2,4)$ is  the set of 2-planes inside a four dimensional space $\C^4$. It is a projective variety  embedded in ${\bP}^5$ and the  embedding is known as the {\it Pl\"{u}cker embedding}. This can be found
explained in many references (see for example Refs. \cite{va,flv}).

Let $\{e_1, e_2,e_3,e_4\}$ be the canonical basis of $\C^4$ and $\pi_0=\rspan\{e_1,e_2\}$ the 2-plane generated by $e_1$ and $e_2$. Then
$$G(2,4)=\rSL(4,\C)/P_0,$$ where the isotropy group $P_0$ is the  subgroup that leaves invariant $\pi_0$:
$$P_0=\left\{\begin{pmatrix}L&M \\
0&R\end{pmatrix}\in \rSL(4,\C)\right\},$$ with
$L, M, R$ being $2\times 2$ matrices, and $L$ and $R$ invertible.

We will call $\rSL(4,\C)$ the
\textit{conformal group} in four dimensions. Indeed, this is the complexified spin
group of $\rSO(2,4)$, easily recognizable as the conformal group of
the Minkowski metric in four dimensions $\rdiag(1,-1,-1,-1)$.
Notice also that we can substitute $\rSL(4,\C)$ by $\rGL(4,\C)$
and the construction will also work.

 To obtain the Pl\"{u}cker map, one starts by considering the vector space $E=\Lambda^2(\C^4)\approx\C^6$ with basis  $\{e_i\wedge
e_j\}_{i<j}$, $\,i,j=1,\dots , 4$. Then, for any plane $\pi=\rspan\{a,b\}$ we have
$$a\wedge b=\sum_{i<j}y_{ij}e_i\wedge e_j.$$
A change of basis $(a',b')=(a,b)u$, where $u\in \rGL(2,\C)$, produces a change
$$a'\wedge b'=\det (u) a\wedge b,$$
so the image $\pi\rightarrow [a\wedge b]$ is well defined into the
projective space $\bP(E)\thickapprox {\bP}^5$. This is called the
Pl\"ucker map and it is not hard to see that it is an embedding. The
quantities $y_{ij}$ are the {\it homogeneous Pl\"ucker coordinates},
and the image of the Pl\"ucker map can be identified with the
solution of the equation  $p\wedge p=0$ for $p\in E$. In coordinates, this reads

\begin{equation}
y_{12} y_{34}+y_{23} y_{14}+y_{31}
y_{24}=0, \qquad (\hbox{ with } y_{ij}=-y_{ji}). \label{klein}
\end{equation}
Equation (\ref{klein}) is called the {\it Pl\"{u}cker relation} or {\it Klein quadric}. In algebraic terms, the fact that the Grassmannian is embedded in the projective space is reflected in the fact that the ideal $\cI_P$  generated by the relation (\ref{klein}) in $\rPol(\C^6)$ is homogeneous.

%\paragraph
{\bf The Poincar\'{e} group plus dilatations and the big cell.}

We give an open cover for the projective variety $G(2,4)$. Let $a$ and $b$ be two linearly independent vectors spanning the 2-plane $\pi$ that can be represented by the matrix
$$\pi=\rspan\{a,b\}=\begin{pmatrix}a_1&b_1\\a_2&b_2\\
a_3&b_3\\a_4&b_4\end{pmatrix}.$$
This matrix has rank 2, so at least one of the
6 minors $y_{ij}=a_ib_j-b_ia_j$, $\,i < j$
is different
from zero. The sets \begin{equation}U_{ij}=\left\{(a,b)\;/\;y_{ij}\neq
0\right\}\label{openset}\end{equation} are open affine sets of $G(2,4)$ and cover it. $U_{12}$ is called the {\it big cell}. By a change of basis (action of $u\in \rGL(2,\C)$), a plane $\pi$ in the big cell can always be represented by the matrix
\be
\pi=\begin{pmatrix}\id\\A\end{pmatrix}, \qquad A=\begin{pmatrix}a_{11}&a_{12}\\a_{21}&a_{22}\end{pmatrix},\label{bigcell}\ee  so $U_{12}\approx M_2(\C)\approx\C^4$, and  $U_{12}$ is again the affine space $\C^4$, represented with the algebra $\C[a_{ij}]$.

Given an element of $\rSL(4,\C)$, since the columns are linearly independent, we can choose the  first two columns to be the vectors $a$ and $b$ representing a 2-plane $\pi$. If the plane is in the big cell, it is easy to see that there is a transformation in $P_0$ that brings the matrix of $\rSL(4,\C)$ to the form
$$\begin{pmatrix}\id_2&0\\A&\id_2\end{pmatrix}.$$

The big cell is left invariant by the subgroup of $ \rSL(4,\C)$ consisting
of the matrices of the form
$$\left\{\begin{pmatrix}L&0\\NL&R\end{pmatrix},\;\; L, R
\;\;\mathrm{invertible},\; \det R\cdot\det L=1\right\}.$$ The bottom
left entry is arbitrary but we have written it like that for convenience. The action on $U_{12}$ is then
\begin{equation}A\mapsto N+RAL^{-1},\label{poincare}\end{equation} so $P$ has the structure of semidirect
product $P=H\ltimes M_2$, where $M_2=\{N\}$ is the set of $2\times 2$ matrices  acting as translations, and
$$H=\left\{\begin{pmatrix}L&0\\0&R\end{pmatrix},\; \; L,R\in
\rGL(2,\C),\;\; \det\!L\cdot \det\!R=1\right\}.$$ The subgroup $H$
is the direct product  $\rSL(2,\C)\times \rSL(2,\C)\times
\C^\times$. But $\rSL(2,\C)\times \rSL(2,\C)$ is the spin group of
$\rSO(4,\C)$, the complexified Lorentz group, and $\C^\times$ acts
as a dilation. $P$ is then the Poincar\'{e} group times dilations. In the basis of the Pauli matrices
\be
\sigma_{0}=\begin{pmatrix}1&0\\0&1\end{pmatrix},\quad
\sigma_{1}=\begin{pmatrix}0&1\\1&0\end{pmatrix},\quad
\sigma_{2}=\begin{pmatrix}0&-i\\i&0\end{pmatrix},\quad
\sigma_{3}=\begin{pmatrix}1&0\\0&-1\end{pmatrix},\label{Pauli}
\ee
an arbitrary matrix $A$ can be written as
$$
A=\begin{pmatrix}a_{11}&a_{12}\\a_{21}&a_{22}\end{pmatrix}=x^{0}\sigma_{0}+x^{1}\sigma_{1}+x^{2}\sigma_{2}+x^{3}\sigma_{3}=
\begin{pmatrix}
x^{0}+x^{3} & x^{1}-ix^{2} \\
x^{1}+ix^{2} & x^{0}-x^{3}
\end{pmatrix},
$$
and $(x^0,x^1,x^2,x^3)$ are the ordinary coordinates of Minkowski space. Moreover,
$$\det A=(x^0)^2-(x^1)^2-(x^2)^2-(x^3)^2.$$ This concludes the interpretation of the big cell as the complexification of the Minkowski space.
$\hbs$

\subsection{The Pl\"ucker embedding for the super Grassmannian
\label{pluckersupergrass}}

We are interested in  the super Grassmannian of $(2|0)$ planes inside
the superspace $\C^{4|1}$. Since we are only concerned with this
particular Grassmannian, we will just denote it as $\Gr=G(2|0;4|1)$.
In the Appendix \ref{supergeometry-app} we give an exhaustive description
of its functor of points, $h_{\Gr}$.  Here, we will use the fact that we can work on local algebras. Then the projective modules are free modules and the description is greatly simplified. On a local algebra $\cA $,  $h_{\Gr}(\cA )$ consists of free submodules
of rank $2|0$ in $\cA ^{4|1}$. One such module can be specified by a couple of independent even vectors, which in the canonical basis $\{e_1,e_2,e_3,e_4,\ep_5\}$ are
\be\pi=\rspan\{a,b\}=\rspan\left\{\begin{pmatrix} a_1&b_1 \\ a_2&b_2\\ a_3&b_3\\ a_4&b_4 \\ \alpha_5&\beta_5 \end{pmatrix}\right\},\label{fopGr}\ee
with $a_i, b_i\in \cA _0$ and $\alpha_5,\beta_5\in \cA _1$.

\bigskip

We want now to work out the expression for the
\textit{Pl\"{u}cker embedding}.
We want to give a natural transformation among the functors
$$
p:h_{\Gr}
\rightarrow h_{\bP(E)}.
$$
where $E$ is the super vector space $E=\wedge^2\C^{4|1}\approx \C^{7|4}$.
We recall that
for an arbitrary super vector space $V$, $$\Lambda^2 V= V\otimes V/ \langle u \otimes v +(-1)^{|u||v|}v \otimes u\rangle, \qquad u,v\in V.$$ Given ,the canonical basis for $\C^{4|1}$ we construct a basis for $E$
\begin{align}&e_1\wedge e_2, e_1\wedge e_3, e_1\wedge e_4, e_2\wedge e_3, e_2\wedge e_4, e_3\wedge e_4, \ep_5\wedge \ep_5,\qquad &\hbox{(even)}\nonumber
\\ &e_1\wedge\ep_5, e_2\wedge\ep_5,e_3\wedge \ep_5,e_4\wedge\ep_5,\qquad &\hbox{(odd)}\label{canonical basis}
\end{align}
As in the super vector space case, if $L$ is a $\cA$-module, we can construct $\wedge^2L $ $$\Lambda^2 L= L\otimes L/ \langle u \otimes v +(-1)^{|u||v|}v \otimes u\rangle, \qquad u,v\in L.$$
If $L \in h_{\Gr}(\cA )$, then  $\wedge^2L \subset \wedge^2\cA ^{4|1}$.
It is clear that if $L$ is a projective $\cA $-module of rank $2|0$, then
$\wedge^2L$ is a projective
$\cA $-module of rank
$1|0$. In other words, it is an element of $h_{\bP(E)}(\cA )$ for
$E=\wedge^2 \C^{4|1}$. Hence we have defined a natural
transformation:
$$
\begin{CD}
h_{\Gr}(\cA ) @>p>> h_{\bP(E)}(\cA ) \\
L @>>> \wedge^2L.
\end{CD}, \qquad \cA \in \salg
$$
Once we have the natural transformation defined, we can restrict to work on local algebras.
For a local algebra $\cA$,
$\Gr(\cA)$ consists of the free $2|0$-modules inside $\cA^{4|1}$.

Let $a,b$ be two even independent vectors in $\cA^{4|1}$. For any
superalgebra $\cA$, they generate a free submodule of $\cA^{4|1}$ of
rank $2|0$. The natural transformation  described above
is as follows.
$$
\begin{CD}
 h_{\Gr}(A) @>p_{\cA}>>  h_{\bP(E)}(\cA) \\
 \rspan_{\cA}\{a,b\} @>>>  \rspan_{\cA}\{a\wedge b\},
\end{CD}
$$
where $E=\Lambda^2(\C^{4|1})$. The map
$p_{\cA}$ is clearly injective. The image $p_{\cA}(h_{\Gr}(\cA))$
is the subset of even elements in $h_{\bP(E)}(\cA)$ decomposable in
terms of two even vectors of $\cA^{4|1}$. We are going to find the necessary
and sufficient conditions for an  even  element $Q\in
h_{\bP(E)}(\cA)$
to be decomposable in terms of even
vectors. In terms of the canonical basis (\ref{canonical basis}) we have
\begin{eqnarray}
&& Q=q+\lambda \wedge \ep_5+a_{55}
\ep_5 \wedge \ep_5 ,\quad
\hbox{with}\nonumber\\
&&q=q_{12}e_1\wedge e_2+\cdots +q_{34}e_3\wedge e_4,\quad
q_{ij}\in \cA_0, \nonumber\\
 &&\lambda=\lambda_1e_1+\cdots +\lambda_4 e_4,\quad \lambda_i\in
 \cA_1.\label{coordinates}
\end{eqnarray}
$Q$ is decomposable if and only if
\begin{eqnarray*}
&&
Q=(r+\xi\ep_5) \wedge (s+\theta\ep_5)\quad \hbox{with}\\
&& r=r_1e_1+\cdots+ r_4e_4,\quad s=s_1e_1+\cdots+ s_4e_4,\quad r_i,
s_i\in \cA_0\quad\xi, \theta\in \cA_1,
\end{eqnarray*} which means
$$
Q=r \wedge s+(\theta r-\xi
 s)\wedge \ep_5+ \xi \theta\ep_5 \wedge \ep_5.$$
This is equivalent to $$q=r \wedge s, \quad \lambda=\theta r-\xi s, \quad
 a_{55}=\xi\theta,
$$
and these  are in turn equivalent to the following:
$$
q \wedge q=0, \qquad q \wedge \lambda=0, \qquad
\lambda\wedge\lambda=2a_{55}q \qquad \lambda a_{55}=0.
$$
Plugging (\ref{coordinates})  we obtain
\begin{align}
&q_{12}q_{34}-q_{13}q_{24}+q_{14}q_{23}=0, &&
\hbox{ (classical Pl\"{u}cker relation)} \nonumber  \\&q_{ij}\lambda_k-q_{ik}\lambda_j+q_{jk}\lambda_i=0,&& 1\leq i<j<k\leq 4\nonumber \\
& \lambda_i \lambda_j=a_{55}q_{ij}&& 1\leq i<j\leq 4\nonumber \\
& \lambda_ia_{55}=0. &&\label{superplucker}
\end{align}
These are the super Pl\"{u}cker relations.
As we shall see in the next section the superalgebra
\be
\cO({\Gr})=k[q_{ij}, \lam_k, a_{55}]/\cI_P, \label{sgr}
\ee
is associated to the supervariety $\Gr$ in the Pl\"{u}cker embedding
described above, where
$\cI_P$ denotes the ideal of the super Pl\"{u}cker relations
(\ref{superplucker}).

\begin{remark} We want to stress that although $\cO({\Gr})$ does
not represent the functor of points $h_{\Gr}$,
(it is not representable as a functor from $\salg$ to $\sets$),
one can recover all the information about the supervariety
from  $\cO({\Gr})$ by the procedure described
in the Appendix \ref{projectivesupergeo}.
(In the notation used there one has that $S= A_{\Gr}$).

$\hbs$\end{remark}

\subsection{The superstraightening algorithm}

We want to take a small digression to explain why the relations
(\ref{superplucker})
generate the ideal $\cI_P$ of all the relations among the coordinates
$q_{ij}$, $\lambda_k$, $a_{55}$. This could be done with a direct calculation
as in \cite{cfl}, however we prefer to justify it by resorting to
the theory of semistandard tableaux and the superstraightening
algorithm. The algorithm has  interest on his own, but we will also
need these notions later, in the quantum setting.

\begin{definition} \label{semistandard}
A \textit{Young diagram} is a finite collection of boxes, or cells,
arranged in left-justified rows, with the row lengths weakly decreasing
(each row has the same or shorter length than its predecessor). Listing
the number of boxes in each row gives a partition of a non-negative
integer $n$, the total number of boxes of the diagram.
A \textit{Young tableau} is obtained by filling in the boxes of the
Young diagram with symbols taken from some alphabet,
which is usually required to be a totally ordered set. Usually the alphabet
consists of the first natural numbers.

Let us assume now that the set of indices is separated into two disjoint subsets:
the {\it even} and the {\it odd} indices. A tableau is called \textit{semistandard}
or \textit{superstandard} the following conditions are satisfied:
\begin{itemize}

\item{ The entries are non
decreasing along each row.}
\item{The rows have no repeated even entries.}
\item{ The entries are non decreasing down each
column.}
\item{ The columns have no repeated odd indices.}\hfill$\blacksquare$
\end{itemize}
\noindent Some authors take the opposite convention, interchanging rows and columns.
\end{definition}

To describe the super Grassmannian we take indices $1$, $2$, $3$, $4$, $5$
with $1$, $2$, $3$, $4$ being even and $5$ being odd. Each of the generators of the super Grassmannian  is associated to a pair of
multi-indices, that we shall write using the letterplace notation \cite{rota}, \cite{brini}. The first pair of indices indicates the rows and the second pair indicates the columns that determine the submatrix whose minor is associated to the generator. Then we have, also in terms of the semistandard tableaux,
$$
\begin{CD}
q_{ij}@>>>(i,j|1,2)@>>> \young(ij)\, ,\, \young(12) \\
\lambda_k @>>>(k,5|1,2)@>>>\young(k5)\, , \,\young(12) \\
a_{55} @>>> (5,5|1,2)@>>>\young(55)\, ,\, \young(12).
\end{CD}
$$
We can suppress the second tableau, which indicates the columns, since it is going to be the same in what follows. A monomial in the generators can be encoded in a tableau
where each line corresponds to the indices of the coordinates, in the order
in which they are written.
For example the monomials $q_{12} q_{34} \lambda_3 a_{55}$,
 $q_{12} q_{23} \lambda_{55}$ correspond
to the tableaux:
$$
{\young(12,34,35,55)} \qquad
{\young(12,23,55)} \qquad
$$
The second one is semistandard, while the first one is not.
As the reader can quickly notice, there is in fact  a relation among
the coordinates whose indices appear in the first tableau:
$\lambda_3a_{55}=0$.

The semistandard tableaux are very important since the monomials associated to them  are
a basis for the superring $\cO(\Gr)$. The next theorem provides us
with a presentation of such ring and with a basis. Its proof is
based on the straightening algorithm in the super setting, which
we are unable to describe here, since it would take us too far
from our purpose. We refer the reader to the beautiful works
\cite{rota}, \cite{brini}, where the full details are discussed.

\begin{theorem}$^{}$\label{pres-grass}

\begin{enumerate}
\item The Grassmannian superring is given in terms of
generators and relations as:
$$
\cO(\Gr)=\C[ q_{ij}, \lambda_{j}, a_{55}]\,/\,\cI_P,
\qquad 1\leq i<j\leq 5 \hbox{  and  } i=j=5
$$
where $\cI_P$ is the two-sided ideal generated by the
super Pl\"{u}cker relations (\ref{superplucker}).

\item The Grassmannian superring is the free superring
over $k$ generated by the monomials in the variables
$q_{ij}$, $\lambda_{j}$, $a_{55}$ whose indeces form a
semistandard tableau.

\end{enumerate}

\hfill$\blacksquare$

\end{theorem}

\begin{remark}
In our special case a semistandard tableaux
means that the indices of the variables forming the monomial
and appearing in the tableau in two consecutive lines are such that
$(i_k,j_k)<(i_{k+1},j_{k+1})$
lexicographically, but the even pairs of type $(1,4)$ and $(2,3)$
(and more generally $(i,j)$, $(k,l)$, with $i<k<l<j$) are not allowed to
appear: such pairs  shall be disposed using the (super)
Pl\"ucker relations, through
the straightening algorithm \cite{rota}.
\hfill$\blacksquare$
\end{remark}

The ring $\cO(\Gr)$ is our starting point for the quantization: we
will quantize the Grassmannian together with its embedding into
$\bP(E)$ by obtaining a quantization of the ring $\cO(\Gr)$. We will obtain the quantized superalgebra as a sub superalgebra of the quantum supergroup $\rSL_q(4|1)$.  In this way, the quantized Grassmannian $\Gr_q$  will carry a natural action of the quantum supergroup $\rSL_q(4|1)$, just as in the classical case.

\subsection{The conformal and Poincar\'{e} supergroups
and the big cell} \label{superpoincare}

We start by describing the natural action of the special linear
supergroup, the conformal supergroup, on the Grassmannian supervariety.

\medskip

The functor of points of the supergroup, in terms of local algebras,  is
\be h_{\rSL(4|1)}(\cA )=\left\{g=\begin{pmatrix}c_{11}&c_{12}&c_{13}&c_{14}&\rho_{15}\\c_{21}&c_{22}&c_{23}&c_{24}&\rho_{25}\\
c_{31}&c_{32}&c_{33}&c_{34}&\rho_{35}\\c_{41}&c_{42}&c_{43}&c_{44}&\rho_{45}\\
\delta_{51}&\delta_{52}&\delta_{53}&\delta_{54}&d_{55}\end{pmatrix},\qquad \mathrm{Ber}(g)=1 \right\},\label{fopSL}\ee where $c_{ij}, d_{55}\in \cA _0$ and $\rho_{i5},\delta_{5i}\in \cA _1$. Ber stands for the
{\it Berezinian} or {\it superdeterminant} of the matrix $g$. We
refer the reader to \cite{cf} Ch. 1 for its definition and main properties.

\begin{remark} The functor $h_{\rSL(4|1)}$ is representable, and it is represented by the algebra
\be\cO(\rSL(4|1))=\C[g_{ij}, g_{55}, \gamma_{i5},\gamma_{5j}]/
(\mathrm{Ber}-1)
\label{algebrasupergroup}\ee
with $i,j=1,\dots, 4$.
Again, latin letters are for even generators and odd letters for odd generators.
If we prefer to use the supergroup $\rGL(4|1)$ we will have instead
$$\cO(\rGL(4|1))=\C[g_{ij}, g_{55}, \gamma_{i5},\gamma_{5j}, X, D]/(X\cdot \det g_{ij}-1, \, D\cdot g_{55}-1),$$ where $X$ and $D$ are even generators.
$\hbs$\end{remark}

We can describe the action of the supergroup  $\rSL(4|1)$ over $\Gr$ as a natural transformation of the functors,
$$
h_{\rSL(4|1)}(\cA )\times h_{\Gr}(\cA )\longrightarrow h_{\Gr}(\cA ),
\qquad \hbox{$\cA$ local,}
$$
which in this language is simply given by the multiplication of matrices (\ref{fopSL}) and (\ref{fopGr}). We leave to the reader as an exercise the
fact that this definition, given for local superalgebras, can
be extended to all superalgebras (see \cite{cf} Ch. 9).
Also, we refer the reader to Appendix \ref{actions} for all
the definitions concerning actions and stabilizers.

As in the even case, we may take $\pi_0=\rspan\{e_1,e_2\}$. The stabilizer
of this point is the upper parabolic sub-supergroup $P_u$, whose functor of points is
\be h_{P_u}(\cA )=\left\{\begin{pmatrix}c_{11}&c_{12}&c_{13}&c_{14}&\rho_{15}\\c_{21}&c_{22}&c_{23}&c_{24}&\rho_{25}\\
0&0&c_{33}&c_{34}&\rho_{35}\\0&0&c_{43}&c_{44}&\rho_{45}\\
0&0&\delta_{53}&\delta_{54}&d_{55}\end{pmatrix}\right\}.\label{upperparabolic}\ee Then, the Grassmannian can be identified with the quotient
$$h_{\Gr}(\cA )=h_{\rSL(4|1)}(\cA )/h_{P_u}(\cA ).$$ The description of homogeneous spaces for super Lie groups is done in detail in \cite{flv}.

These algebras are  commutative  Hopf superalgebras. The comultiplication is given as usual by matrix multiplication, (see for example Refs. \cite{fi5,fi6}, also for the counit and antipode), by organizing the generators in matrix form:
$$\Delta \begin{pmatrix}g_{ij}&\gamma_{i5}\\\gamma_{5j}& g_{55}\end{pmatrix}=\begin{pmatrix}g_{ik}&\gamma_{i5}\\\gamma_{5k}& g_{55}\end{pmatrix}\otimes \begin{pmatrix}g_{kj}&\gamma_{k5}\\\gamma_{5j}& g_{55}\end{pmatrix}.$$
The superalgebra $\cO(\Gr)$ is a sub superalgebra (not a Hopf sub superalgebra) of $\cO(\rSL(4|1))$. It is in fact the superalgebra generated by the corresponding minors, and the Pl\"{u}cker relations are all the relations satisfied by these minors in $\cO(\rSL(4|1))$.

As in the non super case, the super Grassmannian  admits an open cover in terms of affine superspaces. In terms of the functor of points we say that $\Gr$ admits a cover by  `open  affine subfunctors'. This is explained in detail in the Appendix \ref{dg}, and it generalizes the open cover of the non super case given in Section \ref{plucker}.
As we have detailed in the previous section, we shall concentrate
our attention just on local algebras. We will describe $h_{U_{12}}$,
the functor of points of the {\sl big cell} $U_{12}$.
First of all,
we write an element of $h_{\rSL(4|1)}(\cA)$ in blocks as (see (\ref{fopSL}))
$$\begin{pmatrix}C_1&C_2&\rho_1\\C_3&C_4&\rho_2\\\delta_1&\delta_2&d_{55}\end{pmatrix}.$$ Assuming that $\det C_1$ is invertible,
we can bring this matrix, with a transformation of $h_{P_u}(\cA)$, to the form
\be\begin{pmatrix}\id_2&0&0\\A&\id_2&0\\\alpha&0&1\end{pmatrix}.\label{standardform}\ee
The assumption that $\det C_1$ invertible means that the matrix is in
the open set $U_{12}=\Gr \cap V_{12}$, where $V_{12}$ is the
affine open set defined by taking in $\bP(E)$
the coordinate $q_{12} \neq 0$.

So if the column vectors of
$$\begin{pmatrix}C_1\\C_3\\\delta_1\end{pmatrix}$$ with $\det C_1$
invertible represent a $2|0$-module in the big cell, the same module can be represented by a matrix of the form
$$
\begin{pmatrix}\id_2\\A\\\alpha\end{pmatrix},
\qquad A=\begin{pmatrix}a_{11}&a_{12}\\a_{21}&a_{22}\end{pmatrix},\qquad \alpha =(\alpha_1,\alpha_2),
$$
with the entries of $A$ in $\cA_0$ and the entries of $\alpha$ in $\cA_1$. The big cell $U_{12}$ of $\Gr$ is then the affine superspace
associated with the superalgebra
\be\cO(U_{12})=\C[a_{ij},\alpha_j]\approx\C^{4|2}.\label{bigcellsuperalgebra}\ee

We are now interested in the super subgroup of
$h_{\rSL(4|1)}(\cA)$ that preserves the big cell $U_{12}$.
This is the stabilizer functor $\rStab_{U_{12}}$. According to the
Definition \ref{stabilizer} we have
$$\rStab_{U_{12}}(\cA)=\{g \in h_{\rSL(4|1)}(\cA)\;|\;
g \cdot h_{U_{12}}(\cA')\subset h_{U_{12}}(\cA') \hbox{ for all $\cA$-algebras $\cA'$ }\}.$$
As we remarked after Definition \ref{stabilizer},
$\rStab_{U_{12}}$ is not in general representable, and it does not correspond to the
functor of points of a supervariety. Nevertheless,
there exists a subsupergroup  $\rSt_{U_{12}}$ of $\rSL(4|1)$,
whose functor of points is the largest subgroup functor
of $\rStab_{U_{12}}$.

In our case,  $\rSt_{U_{12}}$ is the lower parabolic sub supergroup $P_l$, whose functor
of points is given in suitable coordinates by matrices of the
type
\be
h_{P_l}(\cA)=\left\{\begin{pmatrix}x&0&0\\tx&y&y\eta\\d \tau &d   \xi& d
\end{pmatrix}\right\},\label{lowerparabolic}
\ee
where $x$ and $y$ are even, invertible $2\times 2$ matrices, $t$ is an even,
arbitrary $2\times 2$ matrix, $\eta$ a $2\times 1$ odd matrix and $\tau, \xi$
a $1\times 2$ odd matrix.   $d$ is an invertible even element given by
the superdeterminant equal 1 condition.

Let us see this.  We clearly
have from the definitions
$$
h_{P_l}(\cA) \subset h_{\rSt_{U_{12}}}(\cA).
$$
Since this embedding is functorial, we have an embedding
$P_l \subset \rSt_{U_{12}}$ of the algebraic supergroups. Since any
supergroup over $\C$ is also a supermanifold (see \cite{fi7})
we have that this is also a supermanifold embedding.

Let us look at the superdimension of the stabilizer $\rSt_{U_{12}}$
and $P_l$. The superdimension is well defined since these
are supermanifolds. To compute them one can look at the tangent spaces. We then have
$$
\mathrm{dim} \rSt_{U_{12}} \leq \mathrm{dim}\rSL(4|1) -
\mathrm{dim} U_{12}=4^2+1|2 \cdot 4 - 4|2=13|6,
$$
but
$$
\mathrm{dim} P_l= 2^2+2^2+2^2+1|2+2+2=13|6,
$$
hence $\mathrm{dim} \rSt_{U_{12}}=\mathrm{dim} P_l$.
Now,  the equality $\rSt_{U_{12}}=P_l$ follows from the following
theorem.
\begin{theorem}
Let $M$ and $N$ be two supermanifolds with the
same dimension and such that $|M|=|N|$. If we
have an embedding $M \subset N$ then $M \cong N$.
\end{theorem}
{\sl Proof.} See  in \cite{cf} ch. 5. $\hbs$

\bigskip

The action of the stabilizer supergroup $P_l$ on the big cell $U_{12}$
is as follows,
$$\begin{CD}h_{P_l}(\cA)\times h_{U_{12}}(\cA)@>>>\quad h_{U_{12}}(\cA)\\\\
\left(\begin{pmatrix}x&0&0\\tx&y&y\eta\\d \tau &d \xi& d \end{pmatrix},\begin{pmatrix}\id_2\\A\\\alpha\end{pmatrix}\right)@>>>\begin{pmatrix}\id_2\\A'\\\alpha'\end{pmatrix},
\end{CD}$$
where, using a transformation of $h_{P_u}(\cA)$ to revert the resulting matrix to the standard form (\ref{standardform}), we have
\be\begin{pmatrix}\id_2\\A'\\\alpha'\end{pmatrix}=\begin{pmatrix}\id_2\\y(A+\eta\alpha) x^{-1}+t\\
d (\alpha  +\tau +\xi A) x^{-1}\end{pmatrix}.\label{actionbigcell}\ee
The subgroup with $\xi=0$ is the super Poincar\'{e} group times dilations (compare with Eq. (14) in Ref \cite{flv}). In that case
$$d=\det x\det y.$$

\section{Quantum super Grassmannian}\label{secquantumsupergrass}

Before giving the definition of quantum super Grassmannian, we need
some preliminaries on  quantum supergroups in general.

\subsection{Quantum supergroups}
In this section we follow Manin \cite{ma2}.

\begin{definition} \label{ManinCR}
 The  quantum matrix superalgebra $\rM_q(m|n)$ is defined as
$$
\rM_q(m|n)=_{\mathrm{def}}\C_q \langle x_{ij},\xi_{kl}\rangle/\cI_M
$$
where $\C_q\langle x_{ij},\xi_{kl}\rangle$ denotes the free
superalgebra over $\C_q=\C[q,q^{-1}]$
generated by the even variables
$$x_{ij},\qquad  \hbox{ for }\quad 1 \leq i,j \leq m \quad \hbox{ or } \qquad m+1 \leq i,j \leq m+n.$$
and by the odd variables
\begin{align*}&\xi_{kl}&&  \hbox{for  }\quad 1 \leq k \leq m, \quad m+1 \leq l \leq m+n \\&&&\hbox{or   }\,  m+1 \leq k \leq m+n, \quad 1 \leq l \leq m,\end{align*}
satisfying the relation: $\xi_{kl}^2=0$. $\cI_M$ is an ideal generated by relations that we will describe shortly.

To simplify the notation it is convenient sometimes to have a common notation for even and odd variables
$$
a_{ij}=\begin{cases} x_{ij} & 1 \leq i,j \leq m, \, \hbox{   or   } \quad
                         m+1 \leq i,j \leq m+n \\ \\
              \xi_{ij} &  1 \leq i \leq m,\quad  m+1 \leq j \leq m+n, \, \hbox{   or } \\&
                      m+1 \leq i \leq m+n, \quad 1 \leq j \leq m
\end{cases}
$$
The ideal $\cI_M$ is generated by the relations \cite{ma2}:

\begin{align*}
&a_{ij}a_{il}=(-1)^{\pi(a_{ij})\pi(a_{il})}
q^{(-1)^{p(i)+1}}a_{il}a_{ij}, && \hbox{for  } j < l \\&&& \\
&a_{ij}a_{kj}=(-1)^{\pi(a_{ij})\pi(a_{kj})}
q^{(-1)^{p(j)+1}}a_{kj}a_{ij}, && \hbox{for  } i < k \\ &&&\\
&a_{ij}a_{kl}=(-1)^{\pi(a_{ij})\pi(a_{kl})}a_{kl}a_{ij}, &&  \hbox{for  }
i< k,\;j > l \\&&&\hbox{or } \quad i > k,\; j < l \\&&& \\
&a_{ij}a_{kl}-(-1)^{\pi(a_{ij})\pi(a_{kl})}a_{kl}a_{ij}=(-1)^{\pi(a_{ij})\pi(a_{kj})}(q^{-1}-q)a_{kj}a_{il},&&\\
&&& \hbox{for  }\quad i<k,\;j<l
\end{align*}
where $\pi(a_{ij})=p(i)+p(j)$ mod 2 denotes the parity of $a_{ij}$ (with $p(i)=0$ if $1 \leq i \leq m$ and  $p(i)=1$ otherwise).
$\hbs$
\end{definition}

\begin{remark}In our definition we take $q$ to be an indeterminate. However this
definition makes sense also for any value of $q \in \C^\times$ (usually one
asks away from roots of one). In any case, whenever we make use of
the fact that $q$ needs to be an indeterminate we will say it.$\hbs$\end{remark}

$\rM_q(m|n)$ is a  super bialgebra  with the usual comultiplication and
counit:
\be
\Delta(a_{ij})=\sum a_{ik} \otimes a_{kj},
\qquad
\ep(a_{ij})=\de_{ij}.\label{coalgebra}
\ee
Notice that the comultiplication, which encodes the matrix product law (and then the group law) is not deformed. In fact $\rM(m|n)$ is a bialgebra with the commutative product ($q=1$ in $\cI_M$) and the same comultiplication and counit.
For later use, we will display here the commutative diagram for the  coassociativity axiom satisfied by the comultiplication,
\be\begin{CD}\rM_q(m|n)@>\Delta>>\rM_q(m|n)\otimes \rM_q(m|n)\\
@V\Delta VV @VV\Delta\otimes\id V
\\\rM_q(m|n)\otimes \rM_q(m|n)@>>\id\otimes\Delta>\rM_q(m|n)\otimes \rM_q(m|n)\otimes \rM_q(m|n)\label{coassociativity}
\end{CD}\ee

\medskip

We are ready to define the general linear supergroup. Let
\begin{align*}
D_{1}&=_{\rdef}\sum_{\s \in S_m}(-q)^{-l(\s)}
a_{1\s(1)} \dots a_{m\s(m)} \\ \\
D_{2}&=_{\rdef}\sum_{\s \in S_n}(-q)^{l(\s)}
a_{m+1,m+\s(1)} \dots a_{m+n,m+\s(n)}
\end{align*}
be the quantum determinants of the diagonal blocks.
\begin{definition}
The \textit{quantum general linear supergroup} $\rGL_q(m|n)$ is defined as
$$
\rGL_q(m|n)=_{\mathrm{def}}
\rM_q(m|n)
[ {D_{1}}^{-1},{D_{2}}^{-1}],
$$
where ${D_{1}}^{-1}$ and
${D_{2}}^{-1}$ are even indeterminates such that
$$
\begin{array}{c}
{D_1}D_1^{-1}=1=
{D_1}^{-1}D_{1},\\ \\
{D_{2}}
{D_{2}}^{-1}=1=
{D_{2}}^{-1}
{D_{2}}
\end{array}
$$
The \textit{quantum general linear supergroup} $\rGL_q(m|n)$ is defined as
$$
\rSL_q(m|n)=_{\mathrm{def}}
\rM_q(m|n) \, \big/ \,
\langle {\mathrm Ber_q}-1\rangle,
$$
where ${\mathrm Ber_q}$ is the \textit{quantum Berezinian}
(for its definition and properties refer to \cite{fi6}).
$\hbs$
\end{definition}
$\rGL_q(m|n)$ and $\rSL_q(m|n)$ are Hopf superalgebras.
The comultiplication and the
counit are the same as in $\rM_q(m|n)$.
One must give also the comultiplication on $D_1^{-1}$ and $D_2^{-1}$,
and the antipode $S$. This has been done in detail in Refs. \cite{fi5,fi6}.

The roles of $\rGL_q(m|n)$ and $\rSL_q(m|n)$
can be interchanged in what follows, as it happens for the super non quantum setting.
We then make the choice to use $\rSL_q(m|n)$.

\subsection{Presentation of the quantum super Grassmannian
$\Gr_q$}

Let the notation be as above. We are going to define the quantum super Grassmannian $\Gr_q$ in terms of generators and relations. Then
we will prove that it is a deformation of the algebra $\cO(\Gr)$
given in (\ref{sgr}). We are now ready for the central definition
of our paper, namely the quantum deformation of the super Grassmannian.
\begin{definition}
The \textit{quantum super Grassmannian}
of $2|0$ planes in $4|1$ dimensional superspace is the non commutative
superalgebra $\Gr_q$ generated by the following quantum super minors in
$\rSL_q(4|1)$:
\begin{align*}
&D_{ij}= a_{i1}a_{j2}-q^{-1}a_{i2}a_{j1}, &&1\leq i<j \leq 4, \qquad
 \\
&D_{i5}= a_{i1}a_{52}-q^{-1}a_{i2}a_{51}, &&1\leq i \leq 4\\
 &D_{55}=a_{51}a_{52}.&&
\end{align*}$\hbs$
\end{definition}
\noindent For clarity we write explicitly all the generators:
\begin{align}
&&&D_{12}, \quad D_{13}, \quad D_{14}, \quad D_{23}, \quad
D_{24}, \quad D_{34}, \quad D_{55},&\qquad \hbox{(even)}\nonumber\\ &&& D_{15}, \quad D_{25},
\quad D_{35}, \quad D_{45}&\qquad \hbox{(odd)}.\label{minorsgrassmannian}
\end{align}
Notice that $D_{55}$ is an even nilpotent element.

The parity is easily given by the rule:
$$
|D_{ij}|=(-1)^{|i|+|j|}, \qquad |i|=0 \hbox{ for }
i=1, \dots 4, \hbox{ and }
|i|=1 \hbox{ for }  i=5.
$$
We will see that this subalgebra of $\rSL_q(4|1)$ is generated, as a vector space, by monomials in the above determinants. This fact is not obvious at all, since the
commutation relations of the minors could introduce minors that are not included in
(\ref{minorsgrassmannian}). For $q=1$ we recover the classical
Grassmannian algebra $\cO(\Gr)$ described in detail in Sec.
\ref{pluckersupergrass}.

More precisely, we will find a presentation of $\Gr_q$ in terms of generators $X_{ij}$ (identified with the $D_{ij}$ in (\ref{minorsgrassmannian})) and relations. In order to do so, we first  work out the commutation relations
of the minors. Then, as in the classical setting, there will be additional relations among the generators:
the {\it  quantum super Pl\"{u}cker relations.}

Let us  start with the commutation relations.
In Ref. \cite{fi2}  such commutation relations
are given for the even minors $D_{ij}$, with $1\leq i<j\leq 4$.
As one can readily check, they hold also when just one of
the indices is $5$.
The reason is that
the commutation of one even and one odd variable in the
matrix bialgebra generated by the $a_{ij}$'s is the same as the commutation of two
even variables, and the expression of $D_{kl}$ in terms of $a_{ij}$ is formally
the same.

The commutation relations are as follows:

\begin{itemize}
\item If $i,j,k,l$ are {\sl not}
all distinct we have:
\be
D_{ij}D_{kl}=q^{-1}D_{kl}D_{ij}, \quad (i,j)<(k,l),
\quad 1 \leq  i<j<k<l \leq 5 \label{commrel1}
\ee
where `$<$'     refers to the lexicographic ordering.

\item
If $i,j,k,l$ are  all distinct we have:
\begin{align}
&D_{ij}D_{kl}=q^{-2}D_{kl}D_{ij}, &&  1 \leq i<j<k<l\leq 5,
\nonumber \\ \nonumber\\
&D_{ij}D_{kl}=q^{-2}D_{kl}D_{ij}-(q^{-1}-q)D_{ik}D_{jl},
&& 1 \leq i<k<j< l \leq 5,  \nonumber\\\nonumber \\
&D_{ij}D_{kl}=D_{kl}D_{ij}, && 1 \leq i<k<l<j\leq 5 .\label{commrel2}
\end{align}
\end{itemize}

So the only quantum commutation relations that have to be computed are for
$$
D_{ij}D_{55}, \qquad D_{i5}D_{j5}, \qquad D_{i5}D_{55}.
$$
After some computations one gets (for $1 \leq i <j \leq 4$):
\begin{align}
&D_{ij}D_{55}=q^{-2}D_{55}D_{ij},
\nonumber\\ \nonumber\\
&D_{i5}D_{j5}=-q^{-1}D_{j5}D_{i5}-(q^{-1}-q)D_{ij}D_{55}=
-qD_{j5}D_{i5}
\nonumber\\\nonumber \\
&D_{i5}D_{55}=D_{55}D_{i5}=0.
\label{commrel3}
\end{align}

This concludes the discussion of the commutation relations.
We are ready to tackle the calculation of the quantum super Pl\"{u}cker relations.

Again using the result for the non super setting (Ref. \cite{fi2}) we have
\begin{align}
&D_{12}D_{34}-q^{-1}D_{13}D_{24}+q^{-2}D_{14}D_{23}=0 \nonumber\\ \nonumber\\
&D_{ij}D_{k5}-q^{-1}D_{ik}D_{j5}+q^{-2}D_{i5}D_{jk}=0,
\qquad 1 \leq i<j<k \leq 4 ,\label{quantumsuperplucker1}
\end{align}
which yield a total of 5 relations.
To these we must add the relations, which can be computed
directly:
\be
D_{i5}D_{j5}=qD_{ij}D_{55}, \qquad 1 \leq i<j \leq 4.\label{quantumsuperplucker2}
\ee

The relations (\ref{quantumsuperplucker1},\ref{quantumsuperplucker2}) are the quantum super Pl\"{u}cker relations. If one specifies $q=1$, the superalgebra becomes commutative and the quantum super Pl\"{u}cker relations become the standard ones (\ref{superplucker}).

Now we want to show that the commutation relations together with the quantum super Pl\"{u}cker relations are all the relations among the determinants. We do this in the following proposition.
\begin{proposition}{$ $ }

\begin{enumerate}
\item The quantum Grassmannian superring is given in terms of
generators and relations as:
$$
\Gr_q=\C_q\langle X_{ij}\rangle/I_{\Gr},\qquad 1\leq i<j\leq 5 \hbox{  and  } i=j=5
$$
where $I_{\Gr}$ is the two-sided ideal generated by the
commutations relations (\ref{commrel1},\ref{commrel2},\ref{commrel3}) and the quantum super Pl\"{u}cker relations (\ref{quantumsuperplucker1},\ref{quantumsuperplucker2}) where $D_{ij}$ is substituted by the indeterminates
$X_{ij}$. Moreover $\Gr_q/(q-1) \cong \cO(\Gr)$
(see Section \ref{pluckersupergrass}).

\item The quantum Grassmannian ring is the free ring
over $\C_q$ generated by the monomials in the quantum determinants:
$$
D_{i_1j_1} \cdots D_{i_rj_r}
$$
where $(i_1,j_1), \dots ,(i_r,j_r)$ form a semistandard tableau
(see \ref{semistandard}).
\end{enumerate}

\end{proposition}

{\sl Proof.}
Though the proof is based on the classical result and
it is the same as \cite{fh} we briefly sketch it, since
its importance in our construction.

The generic monomials in the quantum determinants $D_{ij}$ generate
the ring $\Gr_q$ as $\C_q$-module. Using the commutation relations
we can certainly write any monomial as a lexicographically
ordered monomial, then using the quantum Pl\"{u}cker relations
and the straightening algorithm \cite{fu}, we can rewrite any lexicographically
ordered monomial  as a linear combination of standard monomials.

Notice that there are two obstacles to apply the straightening
algorithm to the quantum setting, but they are both easily overcome.
The first is the presence of the coefficients $q$ in the Pl\"{u}cker relations.
This is not a problem, since $q$ is invertible.
The second is the
non commutativity: when we are commuting two quantum determinants
one may argue that the pairs of the forbidden kind, that is
$(i,j)$ and $(k,l)$ with $i<k<l<j$ may arise. However a closer
look to the commutation relations shows that this is never the case.

\medskip

So both (1) and (2) will be done if we can show that the standard
monomials in the $D_{ij}$'s are linearly independent (i.e. there
are no other relations among them apart the $I_{\Gr}$ ones).
Assume there is a relation $R$ among such monomials. Clearly
$R= 0$ mod $(q-1)$ since
there are no relations among the standard classical
monomials, hence $R=(q-1)R'$. Such relation evidently holds
also in the bigger superalgebra $\rSL_q({m|n})$, which is known to be torsion
free. Hence the quantum determinants satisfy $R'$ and repeating
this same argument enough times we obtain a non trivial relation
among the classical monomials, hence the relation $R$ we start
with cannot exist.
$\hbs$

\subsection{$\Gr_q$ as a homogeneous quantum space}\label{grqhs}

We want to prove that the quantum super Grassmannian that we have constructed
admits a coaction of a quantum group on it, namely the
quantum group $\rSL_q(4|1)$, as it happens in the classical setting.

\begin{proposition}
$\Gr_q$ is a quantum homogeneous superspace for the quantum supergroup
$\rSL_q(4|1)$, i. e., there is a coaction on $\Gr_q$ given via the restriction of the
comultiplication on $\rSL_q(4|1)$ (\ref{coalgebra}):
$$
\Delta|_{\Gr_q}: \Gr_q \lra  \rSL_q(4|1) \otimes \Gr_q $$
\end{proposition}

{\sl Proof.} We just have to check that the restriction is well defined. Then, the coaction property
$$\begin{CD}\Gr_q@>\Delta|_{\Gr_q}>>\rSL_q(4|1) \otimes \Gr_q\\
@V\Delta|_{\Gr_q} VV @VV\Delta\otimes\id V
\\\rSL_q(4|1)\otimes _{\Gr_q}@>>\id\otimes\Delta|_{\Gr_q}>\rSL_q(4|1)\otimes \rSL_q(4|1)\otimes \Gr_q
\end{CD}$$
is guaranteed by the coassociativity (\ref{coassociativity}).

So we need to verify that
$$
\D(D_{ij}), \quad \D(D_{i5}),\quad  \Delta(D_{55}) \quad \in \quad  \rSL_q(4|1)\otimes \Gr_q.
$$
Let us denote the generic $2\times 2$ quantum minors as
\be
D_{ij}^{kl}=a_{ik}a_{jl}-q^{-1}a_{il}a_{jk},\label{qdet}
\ee
so in the previous notation $D_{ij}=D_{ij}^{12}$. In the purely even setting, we can use the formula  (see Ref. \cite{fi2})
\be
\D(D_{ij}^{12})=\sum_{1\leq k<l\leq 4} D^{kl}_{ij} \otimes  D^{12}_{kl}.\label{coaction}
\ee In the super case, we can extend the sum to $l=5$. Also for the minors   $D_{i5}$ (the only odd ones) the calculation is the same.
This proves immediately that
$$
\D(D_{ij}),\; \D(D_{i5}) \,\in \, \rSL_q(4|1)\otimes \Gr_q, \; \hbox{ for } 1\leq i<j\leq 4.
$$
For $\D(D_{55})$ it is a straightforward (long) check.
$\hbs$

\section{Quantum deformation of the big cell inside the super Grassmannian} \label{secstarproduct}

We want now to define the analogue, in the quantum setting,  of the superalgebra representing the big cell  of the
 Grassmannian
supermanifold. At the classical level we obtained it in (\ref{bigcellsuperalgebra}).

The superalgebra $\cO(U_{12})\approx \C^{4|2}$ is the superalgebra
corresponding to the {\it chiral Minkowski superspace} (see Sec.
\ref{chiralsuperfields}).
In Section \ref{superpoincare} we wrote the action of the lower
parabolic supergroup $P_l$ (which includes the super Poincar\'{e}
group times dilations)  using the functor of
points (\ref{actionbigcell}). We want  now to translate it into the
coaction language
to make the generalization to the quantum setting.
The first step is to understand the Hopf superalgebra of
the  lower parabolic supergroup $P_l$, as the quotient, by a suitable ideal, of the algebra representing $\rSL(4|1)$
$$\cO(\rSL(4|1))=\C[g_{ij}, g_{55}, \gamma_{i5},\gamma_{5j}]/
({\mathrm Ber}_q-1)
$$
(see Eq. (\ref{algebrasupergroup})).
The generators can be written in matrix form
$$\begin{pmatrix}g_{ij}&\gamma_{i5}\\\gamma_{5j} &g_{55}\end{pmatrix},$$
so to read the comultiplication as the matrix product.
\begin{proposition}\label{parabolic}
Let $\cO(P_l)$ be the superalgebra:
$$
\cO(P_l):=\cO(\rSL(4|1))/\cI
$$
where $\cI$ is the (two-sided) ideal generated by
$$g_{1j},g_{2j}, \quad \hbox{for}\quad j=3,4 \quad \hbox{and}\quad \gamma_{15},\gamma_{25}.$$ This is
the Hopf superalgebra of the lower parabolic subgroup, with
comultiplication naturally inherited by
$\cO(\rSL(4|1))$.
\end{proposition}
{\sl Proof.}
One can readily check that $\cI$ is a Hopf ideal, in other words
$$
\Delta(\cI)=\cI \otimes \cO(\rSL(4|1)) + \cO(\rSL(4|1)) \otimes \cI.
$$
Hence the Hopf superalgebra structure goes to the quotient.
The fact that $\cO(P_l)$ represents the lower parabolic supergroup $P_l$ is also clear.
$\hbs$

\medskip

In matrix form, for $\cA$ local, we have
\be
h_{P_l}(\cA)=\left\{
\begin{pmatrix}
g_{11} & g_{12} & 0 & 0 & 0 \\
g_{21} & g_{22} & 0 & 0 & 0 \\
g_{31} & g_{32} & g_{33} & g_{34} & \ga_{35} \\
g_{41} & g_{42} & g_{43} & g_{44} & \ga_{45} \\
\ga_{51} & \ga_{52} & \ga_{53} & \ga_{54} & g_{55} \\
\end{pmatrix}
\right\} \subset h_{\rSL(m|n)}(\cA).\label{lowparabolic}
\ee
The superalgebra representing the big cell
 is in fact a sub superalgebra (not a Hopf subalgebra) of $\cO(P_l)$.
It is more convenient to make a change of variables for the generators, so
$$
\begin{pmatrix}
g_{11} & g_{12} & 0 & 0 & 0 \\
g_{21} & g_{22} & 0 & 0 & 0 \\
g_{31} & g_{32} & g_{33} & g_{34} & \ga_{35} \\
g_{41} & g_{42} & g_{43} & g_{44} & \ga_{45} \\
\ga_{51} & \ga_{52} & \ga_{53} & \ga_{54} & g_{55} \\
\end{pmatrix}=\begin{pmatrix}
x      & 0 & 0 \\
tx     & y & y\eta\\
\ttau x & d\xi & d \\
\end{pmatrix}.
$$
The notation used here is slightly different to the notation used in
(\ref{lowerparabolic}). We can define
\be d\tau=\ttau x,\label{ttau}\ee
but we will see that having $\ttau$ is essential to describe the bigcell.
With this change of variable we have:
$$
h_{P_l}(\cA)\,=\,\begin{pmatrix}
x      & 0 & 0 \\
tx     & y & y\eta\\
d\tau & d\xi & d \\
\end{pmatrix}.
$$

\begin{proposition}
The Hopf superalgebra $\cO(P_l)$ is generated by the two alternative set
of variables:
\begin{itemize}
\item $x$, $y$, $t$, $\ttau$, $\xi$, $\eta$ and $d$;
\item $x$, $y$, $t$, $\tau$, $\xi$, $\eta$ and $d$.
\end{itemize}
The sub superalgebra of  $\cO(P_l)$ generated by $(t, \ttau)$
coincides with  the big cell superring
$\cO(U_{12})$ as defined in (\ref{bigcellsuperalgebra}).

Moreover, there is a well defined
coaction of $\cO({P_l})$ on $\cO(U_{12})$ induced by the comultiplication
(\ref{comultiplicationPl}),
$$
\begin{CD}
\tilde\Delta:\cO(U_{12})@>>>\cO({P_l}) \otimes \cO(U_{12})
\end{CD}
$$
which explicitly takes the form:
\begin{align*}
\tilde\Delta t_{ij}&=t_{ij}\otimes 1+y_{ia} S(x)_{bj}\otimes t_{ab}+
y_{i}\eta_aS(x)_{bj}\otimes \ttau_{jb},&\nonumber \\ \\
\tilde\Delta\ttau_j=&(d\otimes 1)(\tau_a\otimes 1+
\xi_b\otimes t_{ba}+1\otimes \ttau_a)(S(x)_{aj}\otimes 1).&\nonumber
\end{align*}
(The reader should notice right away that this is the dual to
the expression \ref{actionbigcell}).
\end{proposition}

{\sl Proof.}
First of all it is clear that the two given sets generate the
superring $\cO(P_l)$.
The comultiplication of the new variable introduced
can be computed in terms of the comultiplication of the old variables
\be\Delta\begin{pmatrix}
x      & 0 & 0 \\
tx     & y & y\eta\\
\ttau x & d\xi & d \\
\end{pmatrix}=\begin{pmatrix}
x      & 0 & 0 \\
tx     & y & y\eta\\
\ttau x & d\xi & d \\
\end{pmatrix}\otimes \begin{pmatrix}
x      & 0 & 0 \\
tx     & y & y\eta\\
\ttau x & d\xi & d \\
\end{pmatrix}.\label{comultiplicationPl}\ee
From (\ref{comultiplicationPl}) we have that
\begin{align*}
\Delta x&=x\otimes x,&\\ \\
\Delta (tx)&= tx\otimes x+y\otimes tx+y\eta\otimes\tau x,&\\ \\
\Delta (\ttau x)&=\tau x\otimes x+d\xi\otimes tx+d\otimes\tau x.&
\end{align*}
It is convenient to write this in component form. From (\ref{comultiplicationPl}) we have that
\begin{align*}
\Delta x_{ij}&=x_{ik}\otimes x_{kj},&\\ \\
\Delta (t_{ij}x_{jl})&= (\Delta t_{ij})(\Delta x_{jl})=
t_{ia}x_{ab}\otimes x_{bl}+y_{ia}\otimes t_{ab}x_{bl}+
y_{ia}\eta_a\otimes\ttau_{jb} x_b=
&\\&\quad(t_{ij}\otimes 1+y_{ia} S(x)_{bj}\otimes t_{ab}+
y_{i}\eta_aS(x)_{bj}\otimes \ttau_{jb})(x_{jp}\otimes x_{pl}),&\\ \\
\Delta (\ttau_j x_{jl})&=(\Delta \ttau_j)\otimes(\Delta x_{jl})=
\ttau_j x_{jk}\otimes x_{kl}+d\xi_j\otimes t_{jk}x_{kl}+
d\otimes\ttau_j x_{jl}=&\\
&\quad (\ttau_j\otimes 1+d\xi_bS(x)_{aj}\otimes t_{ba}+
dS(x)_{aj}\otimes \ttau_a)(x_{jk}\otimes x_{kl}),&
\end{align*}
where $S$ is the antipode,
$$S(x)= x^{-1}=\frac 1{\det(x)}\begin{pmatrix} x_{22}&-x_{12}
\\-x_{21}& x_{11}\end{pmatrix}.$$
From these equations we can read the coaction of the group on the big cell,
\begin{align}\tilde\Delta t_{ij}&=t_{ij}\otimes 1+y_{ia} S(x)_{bj}\otimes t_{ab}+y_{i}\eta_aS(x)_{bj}\otimes \ttau_{jb},&\nonumber
\\ \\
\tilde\Delta\ttau_j&=\ttau_j\otimes 1+d\xi_bS(x)_{aj}\otimes t_{ba}+dS(x)_{aj}\otimes \ttau_a.&\label{coaction on bigcell}\end{align}

It is straightforward now to compare (\ref{coaction on bigcell}) with (\ref{actionbigcell}) to realize that the action and the coaction are dual to each other. It is enough to use (\ref{ttau}) but only in the first factor of the tensor product:
\begin{align*}\tilde\Delta\ttau_j=&d\tau_aS(x)_{aj}\otimes 1+d\xi_bS(x)_{aj}\otimes t_{ba}+dS(x)_{aj}\otimes \ttau_a=&\\&(d\otimes 1)(\tau_a\otimes 1+\xi_b\otimes t_{ba}+1\otimes \ttau_a)(S(x)_{aj}\otimes 1),&\end{align*} and we obtain the coaction in the desired form.

$\hbs$

\bigskip

We now turn to the quantum setting. We shall repeat all the
classical arguments, exerting however extreme care, since in all of
our calculations, the Manin commutation relations
(see Definition \ref{ManinCR}) now play
a key role. In order to keep our notation minimal we use the same letters
as in the classical case  to denote the generators of the
quantum big cell and the quantum supergroups.

\begin{proposition}\label{qparabolic}
Let $\cO(P_{l, q})$ be the superalgebra:
$$
\cO(P_{l, q}):=\cO(\rSL_q(4|1))/\cI_q
$$
where $\cI_q$ is the (two-sided) ideal in $\cO(\rSL_q(4|1))$ generated by
\be g_{1j},g_{2j}, \quad \hbox{for}\quad j=3,4 \quad \hbox{and}\quad \gamma_{15},\gamma_{25}.\label{ideal}\ee This is
the Hopf superalgebra of the quantum lower parabolic subgroup, with
comultiplication the one naturally inherited from
$\cO(\rSL_q(4|1))$.
\end{proposition}
{\sl Proof.}
Notice that the comultiplication is the same than in the classical case, so  $\cI_q$ is a Hopf ideal (see Proposition \ref{parabolic})
and the Hopf superalgebra structure goes to the quotient.
$\hbs$

\medskip

\begin{remark}\label{CRideal}The quantum lower parabolic supergroup is generated by the images in the quotient of the generators $g_{ij}$ and $\gamma_{ij}$ that are not listed in (\ref{ideal}). In the quantum case, this is a non trivial fact, because in  the commutation relations among the generators of the ideal may appear generators other than  the ones in (\ref{ideal}), giving then a `bigger' ideal than in the classical case. One can check that this does not happen here (see for example Ref. \cite{cfg}).$\hbs$\end{remark}

 As in the classical case, it is convenient to change coordinates
\be\begin{pmatrix}
g_{11} & g_{12} & 0 & 0 & 0 \\
g_{21} & g_{22} & 0 & 0 & 0 \\
g_{31} & g_{32} & g_{33} & g_{34} & \ga_{35} \\
g_{41} & g_{42} & g_{43} & g_{44} & \ga_{45} \\
\ga_{51} & \ga_{52} & \ga_{53} & \ga_{54} & g_{55} \\
\end{pmatrix}  =\begin{pmatrix}
x      & 0 & 0 \\
tx     & y & y\eta \\
\ttau x & d\xi & d  \\
\end{pmatrix}.\label{qchange}
\ee
Notice that in $\cO(P_{l,q})$ the elements $D_{12}$
and $D_{34}^{34}$ are invertible (these are the quantum determinants (\ref{qdet}). One can compute explicitly the inverse change of variables,
\begin{align*}
&x =\begin{pmatrix}
g_{11} & g_{12}\\ g_{21} & g_{22}  \end{pmatrix},
&&
t =
\begin{pmatrix} -q^{-1}D_{23}D_{12}^{-1} & D_{13}D_{12}^{-1}\\
-q^{-1}D_{24}D_{12}^{-1} & D_{14}D_{12}^{-1}  \end{pmatrix}
\nonumber \\\nonumber \\
&y = \begin{pmatrix}
g_{33} & g_{34}\\ g_{43} & g_{44} \end{pmatrix},&&
d =g_{55}, \\ \\
&\ttau=\begin{pmatrix}g_{55}^{-1}\ga_{51} &
g_{55}^{-1}\ga_{51}  \end{pmatrix}, &
&\xi=\begin{pmatrix}g_{55}^{-1}\ga_{53} &
g_{55}^{-1}\ga_{54}  \end{pmatrix}
\\ \end{align*}
$$\rho =y^{-1} \begin{pmatrix} \ga_{35} \\ \ga_{45}  \end{pmatrix}=
{D_{34}^{34}}^{_1}
\begin{pmatrix} g_{44} & -q^{-1}g_{34}\\ -qg_{43} & g_{33}
 \end{pmatrix}=
\begin{pmatrix} {D^{34}_{34}}^{-1}\ga_{35} \\
{D^{34}_{34}}^{-1}\ga_{45} \end{pmatrix}$$

\noindent It is not hard to see that $\cO(P_{l,q})$ is also generated by
$x, y, d,\eta, \xi$ and $\ttau$.

The {\it quantum Poincar\'{e}  supergroup times dilations} is the quotient of $\cO(P_{l,q})$  by the ideal $\xi=0$. One can also check that it is a Hopf ideal, so the comultiplication goes to the quotient. The quantum  Poincar\'{e} supergroup times dilations is then generated by the images in the quotient of
$x, y, d,\eta$ and $\ttau$. Here, the Remark \ref{CRideal} applies as well. In matrix form, one has
$$\begin{pmatrix}
x      & 0 & 0 \\
tx     & y & y\eta \\
\ttau x & 0 & d  \\
\end{pmatrix}.$$

\begin{definition}We  definethe \textit{quantum big cell}  $\cO_q(U_{12})$
as the subring of $\cO(P_{lq})$
generated by $t$ and $\ttau$.$\hbs$
\end{definition}

We compute now the quantum commutation relations among
the generators of the quantum big cell.

\begin{proposition}
The quantum big cell superring $\cO_q(U_{12})$ has
the following  presentation:
$$
\cO_q(U_{12}) := \C_q \langle t_{ij}, \ttau_{kl} \rangle \, \big/ \, I_U
$$
where $I_U$ is the ideal generated by the relations:
$$
\begin{array}{c}
t_{i1}t_{i2}=q \, t_{i2}t_{i1}, \qquad
t_{3j}t_{4j}=q^{-1} \, t_{4j}t_{3j},
\qquad 1 \leq j \leq 2, \quad 3 \leq i \leq 4
\\ \\
t_{31}t_{42}=t_{42}t_{31}, \qquad
t_{32}t_{41}=t_{41}t_{32}+(q^{-1}-q)t_{42}t_{31},
\\ \\
\ttau_{51}\ttau_{52}=-q^{-1}\ttau_{52}\ttau_{51}, \qquad
t_{ij}\ttau_{5j}=q^{-1}\ttau_{5j}t_{ij}, \qquad 1 \leq j \leq 2 \\ \\
t_{i1}\ttau_{52}= \ttau_{52}t_{i1},
\qquad  t_{i2}\ttau_{51}=\ttau_{51} t_{i2}+(q^{-1}-q)t_{i1}\ttau_{52}.
\end{array}
$$
\end{proposition}
{\sl Proof.} Direct check. $\hbs$

\medskip

As in the classical setting we have the following proposition.

\begin{proposition}
The quantum big cell $\cO_q(U_{12})$ admits a coaction of
$\cO(P_{l,q})$ obtained by
restricting suitably the comultiplication in
$\cO{(P_{l,q})}$. Explicitly (see
\ref{coaction on bigcell})),
\begin{align*}\tilde\Delta t_{ij}&=t_{ij}\otimes 1+y_{ia} S(x)_{bj}\otimes t_{ab}+y_{i}\eta_aS(x)_{bj}\otimes \ttau_{jb},&
\\ \\
\tilde\Delta\ttau_j&=(d\otimes 1)(\tau_a\otimes 1+\xi_b\otimes t_{ba}+1\otimes \ttau_a)(S(x)_{aj}\otimes 1).&
\end{align*}
by choosing as before the set of generators $x$, $y$, $t$, $d$, $\tau$, $\rho$ and $\xi$
for $\cO(P_{l,q})$ and $t$, $\ttau$ for $\cO_q(U_{12})$ with
$d\tau=\ttau x$.
\end{proposition}

{\sl Proof.}
This is so because the comultiplication is the same in the classical
and the quantum group, given essentially by matrix multiplication.
One has to be careful, though,  when expressing the comultiplication
in terms of the new  generators, since the ordering appearing in the
Definition (\ref{qchange}) has to be kept consistently.
$\hbs$

\section{Conclusions}\label{secconclusions}

In this paper we have obtained a quantum chiral conformal superspace as a non commutative superalgebra that admits the action of the superconformal group $\rSL_q(4|1)$. The quantum chiral Minkowski superspace is realized as the big cell inside the quantum conformal superspace, and the quantum super Poincar\'{e} group is properly defined as a sub supergroup of the conformal supergroup that preserves the bigcell.

In particular, we have used the non obvious property of the chiral conformal superspace, the super Grassmannian of $2|0$-planes inside $\C^{4|1}$, of having an embedding in a superprojective space. The quantization has been performed using explicitly this embedding, thus giving implicitly a deformation of the projective superspace perhaps to be compared with the one in Ref. \cite{imt}.

 To obtain a supermanifold
which admits the correct real
form of the Minkowski and conformal superspaces, one has to go to a larger supergeometric object,
namely the {\it flag supermanifold} $F(2|0,2|1,4|1)$
\cite{ma1,kn,flv}. Luckily enough, this  superflag is also projective (which is not true for an arbitrary superflag, contrary to the non super case), so the same method employed here to quantize  the Grassmannian can be used for the superflag. This, however, is work that has its own peculiarities and that we leave for a further paper.

From the physical point of view, this opens the possibility of constructing supersymmetric field theories with chiral superfields on genuinely noncommutative superspaces which carry an undeformed action of the supersymmetry algebra.

\section*{Acknowledgments}

We want to thank V. S. Varadarajan for his help in the
elaboration of the paper and for his hospitality at the
Department of Mathematics at UCLA. We also wish to thank prof. L.
Migliorini for helpful comments.

D. Cervantes wants to thank the Departament de F\'{\i}sica Te\`{o}rica,
Universitat de Val\`{e}ncia for the hospitality
during the elaboration of this work.

M. A. Lled\'{o} wants to thank the Universit\`a degli Studi di Bologna for its hospitality during the realization of this work.

This work has been supported in part by grants FIS2008-06078-C03-02 and FPA2008-03811-E/INFN
of Ministerio de Ciencia e Innovaci\'{o}n (Spain) and ACOMP/2010/213 form Generalitat Valenciana.

\appendix

\section{Supergeometry} \label{supergeometry-app}

In this appendix we want to recall some basic definitions and facts
in supergeometry. For more details see refs. \cite{va, fl2, dm, ma1}.

\subsection{Basic definitions}\label{apbasicdefinitions}
For definiteness, we take the ground field to be  $k=\R,\C$. A {\it
superalgebra} $A$ is a $\Z_2$-graded algebra, $A=A_0 \oplus A_1$.
$A_0$ is an algebra, while $A_1$ is an $A_0$-module. Let  $p(x)$
denote the parity of a homogeneous element $x$,
$$p(x)=0\; \hbox{  if }\; x\in A_0,\qquad p(x)=1 \;\hbox{  if }\; x\in A_1.$$
The superalgebra $A$ is said to be {\it commutative} if for any two
homogeneous elements $x, y$
$$
xy=(-1)^{p(x)p(y)}yx.
$$
The category of commutative superalgebras will be denoted by $\salg$.
We call $A_r=A/I_{\mathrm{odd}}$, with $I_{\mathrm{odd}}$ the (two-sided) ideal generated by
the odd nilpotents the, \textit{reduced algebra} associated with $A$.
Notice that $A_r$ may have even nilpotents, making the terminology a bit
awkward.

\medskip

From now on all superalgebras are assumed to be commutative
unless otherwise specified.

\medskip

\begin{definition}
A {\it superspace} $S=(|S|, \cO_S)$ is a topological space $|S|$
endowed with a sheaf of superalgebras $\cO_S$ such that the stalk at
a point  $x\in |S|$ denoted by $\cO_{S,x}$ is a local superalgebra
for all $x \in |S|$, i. e. it has a unique (two-sided) ideal.
\end{definition}

 \begin{definition} A {\it morphism} $\phi:S \lra T$ of superspaces is given by
$\phi=(|\phi|, \phi^\#)$, where $\phi: |S| \lra |T|$ is a map of
topological spaces and $\phi^\#:\cO_T \lra \phi_*\cO_S$ is
a local sheaf morphism, that is,
$\phi_x^\#(\bm_{|\phi|(x)})=\bm_x$, where $\bm_{|\phi|(x)}$ and
$\bm_{x}$ are the maximal ideals in the stalks $\cO_{T,|\phi|(x)}$
and $\cO_{S,x}$ respectively.
\end{definition}

The most important examples of superspaces are given by
{\sl supermanifolds} and {\sl superschemes}. Let us introduce them.

\begin{example}
The superspace $\R^{p|q}$ is the topological space $\R^p$ endowed with
the following sheaf of superalgebras. For any $U
\subset_\mathrm{open} \R^p$
$$
\cO_{\R^{p|q}}(U)=C^{\infty}(\R^p)(U)\otimes \R[\xi^1, \dots, \xi^q],
$$
where $\R[\xi_1, \dots ,\xi_q]$ is the exterior algebra (or {\it
Grassmann algebra}) generated by the $q$ variables $\xi_1, \dots,
\xi_q$. \hfill$\blacksquare$
\end{example}

\begin{definition}\label{supermanifold}
A  \textit{supermanifold} of dimension $p|q$ is a superspace $M=(|M|, \cO_M)$
which is locally isomorphic to the superspace
$\R^{p|q}$, i. e. for all $x \in |M|$
there exist an open set $V_x \subset |M|$ and  $U \subset \R^{p|q}$
such that:
$$
{\cO_{M}}|_{V_x} \cong {\cO_{\R^{p|q}}}|_U.
$$\hfill$\blacksquare$\end{definition}

Let now $S=(|S|,\cO_S)$ be a superspace and
let $\cO_{S,0}$ and $\cO_{S,1}$ denote the following sheaves:
$$
\cO_{S,0}(U):=\left(\cO_{S}(U)\right)_0,\qquad
\cO_{S,1}(U):=\left(\cO_{S}(U)\right)_1,\qquad U\subset_\mathrm{open} |S|.
$$
Notice that $\cO_{S,0}$ is a sheaf of algebras, while
$\cO_{S,1}$ is a sheaf of  $\cO_{S,0}$-modules. We have the
following.
\begin{definition}A {\it superscheme}
$S$ is a superspace $(|S|, \cO_S)$ such that
$\cO_{S,1}$ is a quasi coherent sheaf of
$\cO_{S,0}$-modules.\hfill$\blacksquare$\end{definition}

Morphisms  of supermanifolds  or of superschemes are just the
morphism of the corresponding superspaces.
As for supermanifolds, also superschemes can be characterized by a
local model. Let us briefly describe it.

\begin{definition} $\uspec A$. \label{spec}

Let $A$ be an object of $\salg$. Since $A_0$ is an algebra, we can
consider the topological space $$\spec(A_0)=\{\hbox{prime ideals }
\bp\subset A_0\}$$  with its structural sheaf $\cO_{A_0}$. The stalk $A_{\bp}$
of the structural sheaf at the prime $\bp\in \spec(A_0)$ is the
localization of $A_0$ at $\bp$.
As for any superalgebra, $A$ is a module over $A_0$, and we have
indeed a sheaf $\widetilde A$ of $\cO_{A_0}$-modules over $\spec A_0$
with stalk $A_{\bp}$, the localization of the $A_0$-module $A$ over
each prime $\bp \in \spec(A_0)$.
$\uspec A=_{\defi}(\spec
A_0,\widetilde A)$ is a superscheme.
For more details concerning  the
construction of the sheaf $\widetilde M$ for a generic $A_0$ module $M$,
see Ref. \cite{ha} II \S 5.
\hfill$\blacksquare$
\end{definition}

It is not hard to see that $\uspec$ is also a functor from
the category of superalgebras to the category of superschemes.
The proof is very similar to the ordinary setting (see \cite{eh} Ch. II)
and can be found in \cite{cf} Ch. 10.

\begin{example}
We define $\bA^{m|n}:=
\uspec k[x_1, \dots, x_m,\xi_1, \dots, \xi_m]$. This superscheme
is called the \textit{affine superspace} of dimension $m|n$.
Its underlying scheme is the affine space $\bA^m$ of dimension $m$.
\end{example}

A superscheme $S$ which is isomorphic to $\uspec A$ for some algebra
$A$ is said to be an {\it affine superscheme}. When
the reduced algebra $A_r=A/I_{\mathrm{odd}}$ is finitely generated and
reduced, that is, it has no nilpotents,
we say that $\uspec A$ is a {\it supervariety}. We have the
following.

\begin{proposition}
A superspace $S$ is a superscheme if and only if it is locally
isomorphic to $\uspec A$ for some superalgebra $A$, i. e. if for
every $x \in |S|$, there exists $U_x \subset|S|$ open such that
$(U_x, \cO_S|_{U_x}) \cong \uspec A$ for some superalgebra  $A$ (that
clearly depends on $U_x$).
\end{proposition}
{\sl Proof.} See Ref. \cite{cf} \S 5.\hfill$\blacksquare$

\medskip
Next we want to introduce the concept of {\it functor
of points} of a superscheme.

\begin{definition}
The {\it functor of points} of a
superscheme $X$ is the representable functor:
$$
\begin{CD}
h_X:\sschemes^o @>>>\sets \\
T @>>>h_X(T)=\Hom(T,X)
\end{CD}
$$
and $h_X(\phi)f=f \circ \phi$ for any morphism $\phi:T \lra S$.
The elements in  $h_X(T)$ are called the \textit{$T$-points} of $X$.
\hfill$\blacksquare$
\end{definition}
(The label `$\,{}^o\,$' means that we are taking the opposite category.)

The following facts detailed in
the next observation are not difficult to prove (See  Ref. \cite{cf}).
They will be important in the sequel.

\begin{observation} \label{facts}${}^{}$

\noindent 1. The functor of points of a superscheme is determined by
its restriction to the category of affine superschemes.

\medskip

 \noindent 2. The category of affine superschemes is equivalent
to the category of affine superalgebras, denoted by $\salg$,
hence the functor of
points of a superscheme can be equivalently
defined also as a functor from
$\salg$ to $\sets$. In other words given a superscheme $X$, we can
equivalently define its functor of points as:
$$
\begin{CD}
h_X:\salg @>>>\sets \\
A @>>>h_X(A)=\Hom(\uspec A, X)
\end{CD}
$$
and $h_X(\phi)f=f \circ \uspec \phi$ for any morphism $\phi:A \lra B$.

The elements in  $h_X(A)$ are called the \textit{$A$-points} of $X$.

When $X$ is itself an affine superscheme, the functor of points can
be written as follows:
$$
\begin{CD}
h_X:\salg @>>>\sets \\
A @>>>h_X(A)=\Hom(\cO(X), A)
\end{CD}
$$
and $h_X(\phi)f=\phi \circ f$ for any morphism $\phi:A \lra B$.

Notice that unless $X$ is an affine superscheme, the functor
$h_X:\salg \lra \sets$ will not be representable.

\medskip

\noindent 3. The functor of points of a superscheme seen as $F:\salg
\lra \sets$ is a local functor i. e. it has the sheaf property.
In other words let $A \in \salg$ and
$(f_i,\; i \in I)=(1)=A$. Let
$\phi_{k}:A \lra A_{f_k}$ be the natural map of the algebra $A$ to its
localization and also $\phi_{kl}:A_{f_k} \lra A_{f_kf_l}$.
Then, for a family
$\al_i\in F(A_{f_i})$, such that $F(\phi_{ij})(\al_i)=F(\phi_{ji})(\al_j)$
there exists $\al\in F(A)$ such that
$F(\phi_i)(\al)=\al_i$.

\medskip

\noindent 4. (Yoneda) Given superschemes $S$ and $T$, the natural
transformations $h_S \lra h_T$ are in one-to-one correspondence with
the superscheme morphisms $S \lra T$. Consequently two superschemes are
isomorphic if and only if their functor of points are isomorphic.

\hfill$\blacksquare$
\end{observation}

\subsection{Projective supergeometry} \label{projectivesupergeo}

We want now to consider projective superschemes and supervarieties.
Recall that to an ordinary graded algebra $R=\oplus_{i\geq 0} R^i$, we can
associate the topological space
$$
\Proj R =\{\hbox{relevant homogeneous prime ideals } \bp\subset R\}
$$ (An ideal is {\it relevant} if it does not include $R^+=\oplus_{i>0} R^i$, see Ref. \cite{ha} pg 116).
 The structural sheaf of  $\Proj R $, denoted by  $\cPO_R$
has stalk at $\bp$:
$$
R_\bp=\left\{\frac f g  \; | \;  f\in R,\,
g\hbox{ homogeneous } \in R-\bp\right\}.
$$
We want to generalize this construction to superalgebras.

\medskip

{\bf Notation}. We shall use the lower indices to indicate the
$\Z_2$-gradation, while the upper indices will indicate the
$\Z$-gradation. When we say `graded' we shall always mean
$\Z$-graded, while for the $\Z_2$-graded objects we shall use the
word `super'. $\hbs$

\medskip

Let us consider a graded superalgebra $A=\oplus_{i\geq 0}A^i$, with
$\pi^i:A\rightarrow A^i$ the natural projection. We always
assume that this $\Z$-grading is compatible with the $\Z_2$-grading $A=A_0+A_1$
i. e.
$$
\pi^i(A_0)=A_0\cap A^i.
$$
Then we have that $A_0$,  the even part of $A$, is a graded algebra,
$$
A_0=\oplus_{i\geq 0}A_0^i,\qquad A_0^i=A_0\cap A^i,
$$
and $A$ is a graded $A_0$-module.

\begin{definition} $\uproj A$.

Let $A$ be a graded superalgebra. Similarly to
Definition \ref{spec} we consider the sheaf of graded
superalgebras $\widetilde A$ on the topological space
$\Proj A_0$ with stalk at $\bp\in \Proj
A_0$
$$
A_\bp=\left\{\frac f g\quad  |\quad  f, g\in A,\,
g\hbox{ homogeneous } \in A_0-\bp\right\}.
$$
One can check that
 $(\Proj A_0, \widetilde A)$ is a
superscheme,   and we will denote it with $\uproj A$ (see \cite{ha}
Ch. II \S 5 for more details). \hfill$\blacksquare$
\end{definition}

Let us see two important examples: the projective superspace
and the projective supervarieties.

\begin{example}
\label{projsup} {\it Projective superspace.}

Consider the graded superalgebra $S=\C[x_0 \dots x_{m}, \xi_1 \dots
\xi_n]$. We want to describe $\uproj S$ explicitly as a
superscheme.

For each $r$, $0 \leq r \leq  m$, we consider the graded superalgebra
$$
S[r]=\C[x_0,\dots ,x_m, \xi_1,\dots, \xi_n][x_r^{-1}],\qquad
\deg(x_r^{-1})=-1.
$$
 The subalgebra $S[r]^0\subset S[r]$ of degree 0 is
\begin{equation}
S[r]^0\approx \C[u_0,\dots,\hat u_r ,\cdots, u_m,\eta_1,\dots
\eta_n],\qquad u_s=\frac {x_s}{x_r},\; \eta_\alpha=\frac
{\xi_\alpha}{x_r}, \label{homogeneous}
\end{equation}
(the label `$\;\,\hat{}\;\,$' means that this generator is omitted). Let
$|U_r|$ be the set of homogeneous prime ideals of $S[r]_0$. Because
of (\ref{homogeneous}), a homogeneous prime ideal of $(S_r)_0$
corresponds to a prime ideal of the subalgebra of degree 0,
$S[r]_0^0$ in $S[r]^0=S[r]_0^0 \oplus S[r]_0^1$,
so $|U_r|=\spec S[r]_0^0$. We denote  by $U_r$ the
superscheme
$$
U_r=\uspec S[r]^0\approx(\spec S[r]_0^0,
\widetilde{S[r]^0}), \quad
\cO_{U_r}=_{\defi} \widetilde{S[r]^0}=
\cO_{\bP^{m}}|_{|U_r|} \otimes
\C[\xi_1, \dots ,\xi_n].
$$
where $\bP^{m}$ is the classical  projective space of dimension $m$ and
$\cO_{\bP^{m}}$ its structural sheaf.
 We
have  then found affine open subsuperschemes $U_r=\uspec S[r]_0^0$ whose
topological spaces cover $\Proj S_0$, and
$$
\cO_{U_r}|_{|U_r| \cap |U_s|}=
\cO_{U_s}|_{|U_r| \cap |U_s|}.
$$
We conclude that there exists a unique sheaf on $\Proj(S)_0$ that we
denote as  $\cO_{\bP^{m|n}}$, whose restriction to $|U_i|$ is
$\cO_{U_i}$ and $\uproj S =(\Proj S_0, \cO_{\bP^{m|n}})$.

We will call the superscheme  $\uproj S =\bP^{m|n}$ the {\it super
projective space} of dimension $m|n$.

\hfill$\blacksquare$ \end{example}

More in general
if $E$ is a super vector space and $\Sym(E)$ the algebra of the
polynomial functions on $E$, we will denote with $\bP(E)$ the
superscheme $\uproj \left(\Sym(E)\right)$.

\hfill$\blacksquare$

\begin{example}{\it  Projective supervarieties.} \label{projsubvar}

Let $I\subset S=\C[x_1 \dots x_m, \xi_1 \dots \xi_n]$ be a
homogeneous ideal; then $S/I$ is also a graded superalgebra and we
can consider the superscheme $\uproj (S/I)$. The natural
surjective morphism
$p:S\rightarrow S/I$ defines a morphism of superschemes
$$
\hat p:\uproj (S/I)\rightarrow \uproj S.
$$
Let us see this more in detail. We clearly have a morphism on the
underlying topological spaces as it happens in the classical case.
In order to build a sheaf morphism, we look at an affine cover of $\Proj
(S/I)_0$. We define the topological space:
\begin{eqnarray*}|V_i|=&&
\Proj\left(\frac {\C[x_0,\dots x_m,\xi_1\dots
\xi_n]}{I}[x_i^{-1}]\right)_0=\\
&&=\spec \left(\frac
{\C[x_0,\dots x_m,\xi_1\dots
\xi_n]}{I}[x_i^{-1}]\right)^0_0= \\
&&=\spec \left(\frac {\C[u_0,\dots,\hat u_i,\dots u_m,\eta_1\dots
\eta_n]}{I_{\mathrm{loc}}}\right),
\end{eqnarray*}
where $I_{\mathrm{loc}}=(I[i])^0_0 \subset \C[u_0,\dots,\hat
u_i,\dots u_m,\eta_1\dots \eta_n]$
in our previous notation. We define the affine superscheme
$$
V_i= \uspec \left(\frac {\C[u_0,\dots,\hat u_i,\dots u_m,\eta_1\dots
\eta_n]}{I_{\mathrm{loc}}}\right)
$$

One can check that the supersheaves $\cO_{V_i}$ are such that
$\cO_{V_i}|_{|V_i| \cap |V_j|}=\cO_{V_j}|_{|V_i| \cap |V_j|}$. Hence
as before there exists a supersheaf $\uproj(S/I)$
on $\Proj (S/I)_0$ whose restriction to $|V_i|$ is $\cO_{V_i}$. The
natural morphisms from $V_i$ to $\bP^{m|n}$ i. e. the maps corresponding
to the superalgebras morphisms:
$$
S[x_i^{-1}] \lra (S/I)[x_i^{-1}]
$$
glue accordingly to give the morphism
$\hat p:\uproj (S/I)\rightarrow \uproj S$.
\hfill$\blacksquare$
\end{example}

\subsection{ The functor of points of projective superspace.} We now want to
understand the functor of points of the projective superspace
$\bP^{m|n}$ and of its subvarieties. The situation is essentially
the same as in the classical case. We will briefly sketch it in the
super setting. (For more details of the classical case see Ref.
\cite{eh} pg 111).

\medskip

In Example \ref{projsup} we have given explicitly the projective
superscheme $\bP^{m|n}$ and an open affine covering. We
want to give now explicitly its functor of points, namely
$$
h_{\bP^{m|n}}(\uspec A)= \Hom(\uspec A,\bP^{m|n}), \quad A \in \salg.
$$
Using the
locality property (see 3 of Facts \ref{facts}), one can
prove
that a morphism $\psi: \uspec A \lra \bP^{m|n}$  is determined by a
family of morphisms $\psi_i:\uspec A_{f_i} \lra \bP^{m|n}$, where
$(f_i,\; i \in I)=(1)=A$, provided they induce the
same map on intersections (see Proposition III-39 in Ref.
\cite{eh}). If we denote by $\phi'_{kl}:\uspec A_{f_kf_l} \lra
\uspec A_{f_k}$, the maps of superschemes induced by
$\phi_{kl}:A_{f_k} \lra \uspec A_{f_kf_l}$ ,
then we have that
$$\psi_{i}\circ \phi'_{ij}=\psi_{j}\circ\phi'_{ji}.$$

We want now to give a description of the morphisms $\uspec A \lra
\bP^{m|n}$ when $A$ is a local superalgebra (Proposition III-36 in Ref.
\cite{eh}).

\begin{proposition}\label{localalgebra}
Let $A$ be a local algebra. Then the morphisms $\uspec A \lra
\bP^{m|n}$ are in one to one correspondence with the set of
$(m+1|n)$-tuples $(a_0 \dots a_m, \al_1 \dots \al_n) \in A^{m+1|n}$
with at least one $a_i$ invertible,  modulo  multiplication by an
invertible element in $A$.
\end{proposition}

{\sl Proof.} The proof is the same as the classical one,
we briefly sketch it. Consider the element
 $(a_0 \dots a_m, \al_1 \dots \al_n) \in A^{m+1|n}$
with $a_i$ ($i$ fixed) a unit.
We want to write the corresponding morphism of superschemes $\uspec A \lra
\bP^{m|n}$. We can write immediately
the morphism of superalgebras
$$
\begin{array}{cccc}
\varphi:& k[U_i]=k[u_1 \dots \hat u_i \dots u_m, \eta_1 \dots \eta_n]
& \lra & A \\
& u_j & \mapsto &  {a_j}/{a_i} \\
& \eta_k & \mapsto &  {\al_k}/{a_i} \\
\end{array}
$$
This is well defined with respect to the equivalence relation.
Moreover it defines a sheaf morphism $\varphi^*:\cO_{U_i} \lra \cO_A$,
with $U_i=\uspec k[U_i]$.
One can check that $\varphi_i^*|_{|U_i| \cap |U_j|}=
\varphi_j^*|_{|U_i| \cap |U_j|}$
hence the $\varphi^*_i$ define a morphism $\varphi: \cO_{\bP^{m|n}}\lra \cO_A$.

Conversely, assume that we have a morphism $\psi:\uspec A \lra
\bP^{m|n}$. The topological
map sends the maximal ideal  $\fm$ of $A_0$ into some $|U_i|$,
$|\psi|(\fm) \in |U_i|$. We claim that all $\spec A_0$ is mapped
inside $|U_i|$.
In fact $|\psi|^{-1}(|U_i|)$ is an open set containing
$\fm$ hence necessarily the whole $\spec A_0$.
Hence $|\psi|(\spec A_0)=|U_i|$. The
morphism $\psi$ is then given by a superalgebra map:
$$
\begin{CD}
k[U_i]=k[x_0/x_i, \dots, x_m/x_i, \xi_1/x_i, \dots, \xi_n/x_i] @>>> A \\ \\
{x_k}/{x_i} ,\; {\xi_\mu}/{x_i}@>>> b_k ,\; \beta_l
\end{CD}
$$
So this map determines: $(b_0, \dots, b_i=1, \dots b_m,\beta_1 ,\dots,
\beta_n)$. One can check that this is well defined. In fact if
$|\psi|(\fm)$ is also in $|U_j|$ we have that it corresponds to a
different $m+1|n$-uple which is a multiple of  the previous one.
\hfill$\blacksquare$

\medskip

The next observation characterizes the functor of points
of subvarieties of the projective superspace, that is \textit{projective
supervarieties} as they are defined in Example \ref{projsubvar}.

\begin{observation}
Let $X \subset  \bP^{m|n}$ be a projective supervariety defined
by the homogeneous ideal $I=(f_1 \dots f_r)$.
As we saw in Proposition \ref{localalgebra}, for $A$
a local superalgebra,
the $A$-points of $\bP^{m|n}$ are  the
$(m+1|n)$-tuples $(a_0 \dots a_m, \al_1 \dots \al_n) \in A^{m+1|n}$
with at least one $a_i$ invertible.
We have that such an $(m+1|n)$-tuple corresponds to an $A$-point of
$X$ if and only if it is a zero of all the polynomials
$f_1 \dots f_r$. We leave this to the reader as an exercise.
For more details on the ordinary setting, see \cite{eh} Ch. III.$\hbs$
\end{observation}

\subsection{\label{dg}The functor of points of the Grassmannian superscheme}

In this section we want to construct the functor
of points of the Grassmannian superscheme
(for a general and more detailed treatment of functor
of points and Grassmannian see \cite{cf} \S 3).

\medskip

We want to stress the importance
of the functorial treatment since it is more general
and it gives geometric intuition to the problem. In fact
it allows to recover the description of Grassmannian superscheme
as the set of submodules of rank $r|s$ inside some
free $m|n$-module.
Let us see more in detail this construction.

\medskip

Consider the functor $\Gr: \salg \lra \sets$, with $\Gr(A)$ the
set of projective quotients of rank $r|s$ of $A^{m|n}$, that is,
\begin{equation}
\Gr(A)=\{\alpha:A^{m|n}\rightarrow L,
 \hbox{ with $L$ a projective A-module of rank }
 r|s\} \label{gra}
\end{equation}
where $\al \cong \al'$ if and only if they have the same kernel.

This functor can also be equivalently defined on the objects
as the set:
\begin{eqnarray*}
\Gr(A)&=
\{\hbox{projective submodules $L$ of
$A^{m|n}$ of rank $r|s$}\}.
\end{eqnarray*}
To complete the definition of functor of points
we need to specify $\Gr$ on morphisms $\psi: A \lra B$.

Given  a morphism $\psi:A\rightarrow B$ on $\salg$, we can give to
$B$ the structure of $A$-module by setting
$$
a\cdot b=\psi(a)b,\qquad a\in A,\;b\in B.
$$
Also, given an  $A$-module
$L$, we can construct the $B$-module $L\otimes_AB$. So given
$\psi$ and the element of $\Gr(A)$, $f:A^{m|n}\rightarrow L$, we
have an element of $\Gr(B)$,
$$
\Gr(\psi)(f):B^{m|n}=A^{m|n}\otimes_A B\rightarrow L\otimes_AB.
$$

We want to briefly motivate why
$\Gr$ is the functor of points of a superscheme,
sending the reader to \cite{cf} and \cite{dg} for more details
and the complete proof. We
use the following result \cite{cf} \S 3:

\begin{theorem}\label{theoremfunctors} A functor $F: \salg \lra \sets$
is the functor of points of a superscheme $X$ if and only if
it has the sheaf property and it admits a cover by
open affine subfunctors.
\hfill$\blacksquare$
\end{theorem}

We will start by showing that $\Gr$ admits a cover by open affine
subfunctors. Consider the multiindex $I=(i_1, \dots, i_r| \mu_1,
\dots, \mu_s)$ and the map $\phi_I: A^{r|s} \lra A^{m|n}$ defined
by
\begin{eqnarray*}\phi_I(x_1,\dots ,x_r|\xi_1, \dots ,\xi_s)=&&m|n-\hbox{tuple
with }\\&&x_1, \dots x_r \hbox{ occupying the position } i_1,
\dots, i_r,\\&&
 \xi_1, \dots \xi_s \hbox{ occupying the position } \mu_1,
\dots, \mu_s \\&&\hbox{and the other positions are occupied by
zero}.\end{eqnarray*}
 (For example, let $m=n=2$ and $r=s=1$. Then
$\phi_{1|2}(x,\xi)=(x,0|0, \xi)$).

If is possible to define subfunctors $v_I$ of $\Gr$ that
on local superalgebras look like:
\begin{equation*}v_I(A)=\{\al: A^{m|n} \lra L\;/\;\al \circ \phi_I
\hbox{ is invertible}\}.
\label{subfunctors}\end{equation*}
It turns out that the $v_I$'s are representable and they
correspond to the affine superspace $\bA^{r|s \times m|n}$, moreover
they cover $\Gr$.

\medskip

As for the sheaf property of $\Gr$, notice that
classically we have a functorial
equivalence between the categories of projective finitely generated
$A$-modules and coherent sheaves on $\uspec A$ which are locally
of constant rank (see \cite{ha} pg 111).
One can check that this classical equivalence
translates to the super category, hence our functor
can be identified with:
$$
\Gr(A) \cong \{\cF \subset \cO_A^{m|n}\,/\, \cF \hbox{ is a subsheaf,
of locally costant rank } r|s \}
$$
where $\cO_A^{m|n}=k^{m|n} \otimes \cO_A$.

\medskip

By its very definition this functor is local. Hence we have proven
that $\Gr$ is the functor of points of a superscheme i. e.
$G=h_S$.

\subsection{Actions of Supergroups} \label{actions}

In this section we briefly summarize the definition and main properties of
the stabilizer functor for the action of a supergroup on a superscheme.
Since a supervariety is in particular a superscheme,
all  our statements hold if we replace the word `superscheme' with
`supervariety'.

\begin{definition} \label{action}
We say that a supergroup $G$ acts on a superscheme $X$,
if we have a natural transformation
between their functor of points:
$$
\begin{array}{ccc}
h_G(A) \times h_X(A) & \lra & h_X(A) \\ \\
g, x & \mapsto & g \cdot x
\end{array}
$$
satisfying the usual axioms for an action:

\begin{itemize}

\item $g \cdot (h \cdot x)=gh \cdot x$, for all $g,h \in h_G(A)$,
$x \in h_X(A)$;

\item $1_G \cdot x = x$, for all $x \in h_X(A)$.

\end{itemize}

\end{definition}

Once a supergroup action is defined, we can talk about the
stabilizer of a subfunctor of $h_X$. We are not assuming the subfunctor
to be representable or in general to be a sub-superscheme of $X$.
In the following we are inspired by \cite{ja} Ch. 1, \S 2.

\begin{definition} \label{stabilizer}
Let $Y$ be a subfunctor of $h_X$.
We call $Stab_Y$ the \textit{stabilizer}
of $Y$ the following supergroup subfunctor of $h_G$:
$$
Stab_Y(A):=\left\{g \in h_G(A) \quad \big| \quad g \cdot Y(A')=Y(A')
\quad \forall A-\hbox{algebras} \quad A'\right\}
$$
\end{definition}

Notice that by this definition there is no guarantee that $Stab_Y$ is
the functor of points of a superscheme
and this even in the case in which $Y$ is a subscheme of $X$.
In fact, as it happens already in the ordinary setting,
we have examples for $Y$ an open subscheme of $X$,
for which $Stab_Y$ is not representable. Consider for example the
natural action of the multiplicative group
$G=\bA^1 \setminus \{0\}$ on $X=\bA^1$
by multiplication and let $Y=h_{\bA^1 \setminus \{0\}}$. This is
the functor of points of an
open subscheme of $X$.
Geometrically it is clear that the only point stabilizing the open subset
$\bA^1 \setminus \{0\}$
is the identity. However if one computes the stabilizer of $Y$ one finds that
$ Stab_Y(A)$ consists of the elements $1+n$, with $n$ a
nilpotent element of $A$ and one can prove this
functor is not representable.

\medskip

Despite this complication, we however have some positive answer to
the question whether $Stab_Y(A)$ is
the functor of points of a superscheme (see \cite{cf}
Ch. 6 for more details and the proof of this statement).

\begin{proposition}
Let the notation be as above and assume $Y$ is
the functor of points of a closed subscheme of $X$.
Then $Stab_Y$ is representable and it is a closed subgroup of the
supergroup $G$.
\end{proposition}

\end{document}